\documentclass[a4paper,11pt,twoside]{article}
\usepackage{pst-fill,pst-grad}
\usepackage{textcomp}
\usepackage[titletoc]{appendix}
\usepackage{titlesec}
\usepackage[english]{babel}
\usepackage[utf8]{inputenc} 
\usepackage{graphicx}
\usepackage{amsmath}
\usepackage{float}
\usepackage{fancyhdr}
\usepackage[matrix,arrow,curve]{xy}
\usepackage{pstricks} 
\usepackage{amsmath,amsfonts,verbatim,afterpage,theorem,euscript,mathrsfs,amssymb}
\usepackage{amsfonts}
\usepackage{amssymb}
\usepackage{array}
\usepackage{dsfont}
\usepackage{hyperref}
\usepackage{authblk}

\newtheorem{Definition}{Definition}[section]
\newtheorem{Proposition}{Proposition}[section]

\newtheorem{Lemme}{Lemma}[section]
\newtheorem{Theoreme}{Theorem}

\newtheorem{Remarque}{Remark}[section]

\def \fe{\vec{f}} 
\def \vu{\vec{u}}
\def \vphi{\vec{\varphi}}
\def \P{\mathbb{P}}

\def \Rt{\mathbb{R}^{3}}
\def \Zt{\mathbb{Z}^3}
\def \finpv{\hfill $\blacksquare$  \newline }
\def \pv{{\bf{Proof.}}~}

\def \ds{\displaystyle}

\title{\bf On the Kolmogorov dissipation law in a  damped Navier-Stokes equation}
\author[1]{\small Diego Chamorro\footnote{diego.chamorro@univ-evry.fr}}
\author[2]{\small Oscar Jarr\'in\footnote{\emph{Corresponding author}:\,\,or.jarrin@uta.edu.ec}}
\author[1]{\small Pierre-Gilles Lemari\'e-Rieusset\footnote{plemarie@univ-evry.fr}}

\affil[1]{\scriptsize LaMME, Univ Evry, CNRS, Université Paris-Saclay, 91025, Evry, France.}
\affil[2]{\scriptsize Dirección de investigación y desarrollo (DIDE).
Universidad Técnica de Ambato, 
Avenida de los Chasquis, 180207, Ambato, Ecuador.}

\begin{document}
\maketitle

%%%%%%%%%%%%%%%%%%%%%%%%%%%%%%%%%%%%%%
\begin{scriptsize}
\abstract{We consider here the Navier-Stokes equations in $\mathbb{R}^{3}$ with a stationary, divergence-free external force and with an additional damping term that depends on two parameters. We first study the well-posedness of weak solutions for these equations and then, for a particular set of the damping parameters, we will obtain an upper and lower control for the energy dissipation rate $\varepsilon$ according to the Kolmogorov K41 theory. However, although the behavior of weak solutions corresponds to the K41 theory, we will show that in some specific cases the damping term introduced in the Navier-Stokes equations could annihilate the turbulence even though the  Grashof number (which are equivalent to the Reynolds number) are large. }\\[3mm]
\end{scriptsize}
%%%%%%%%%%%%%%%%%%%%%%%%%%%%%%%%%%%%%%
\section{Introduction} 
In this article we study the Kolmogorov dissipation law in the deterministic framework of the following damped  Navier-Stokes equations
\begin{equation}\label{damped_alpha_kappa_N-S}
\left\lbrace \begin{array}{lc}\vspace{3mm}
\partial_t\vu=\nu\Delta \vu-\P((\vu\cdot \vec{\nabla}) \vu)+\fe-\alpha P_\kappa(\vu), \quad div(\vu)=0,\quad \nu>0,  &\\
\vu(0,\cdot)=\vu_0,&\\
\end{array}\right.
\end{equation}
where $\vu:[0,+\infty[\times \Rt\longrightarrow \Rt$ is the velocity of the fluid, $\P$ is the Leray projector given by  $\P(\vphi)=\vphi-\vec{\nabla} \frac{1}{\Delta}(\vec{\nabla} \cdot \vphi)$, $\nu >0$ is the fluid viscosity parameter, $\fe\in L^2(\Rt)\cap \dot{H}^{-1}(\Rt)$ is a divergence-free, time-independent external force, $\vu_0\in L^2(\Rt)$  is a divergence-free initial data and for the parameters $\alpha>0$ and $\kappa>0$ the term $-\alpha P_\kappa(\vu)$ is a frequency truncation operator defined by formula (\ref{op_P}) below. When $\alpha=0$ we will refer to the problem (\ref{damped_alpha_kappa_N-S}) as the \emph{classical Navier-Stokes equations}. \\

Let us start by explaining the phenomenological idea behind the Kolmogorov dissipation law which is known as the \emph{energy cascade model}. This model explains that the kinetic energy is introduced in the fluid by the effect of the external force $\fe$  at a length scale of order $\ell_0>0$ which is called the \emph{energy input scale} and, as we can control the external force this quantity $\ell_0$ will always be given. Then, in the turbulent  setting of large Reynolds number $Re$ (see  expression (\ref{Reynolds}) below) the energy dissipation mechanism (which is due to the viscosity forces of the fluid)  is not effective and this  energy introduced (at the input scale $\ell_0$) is transferred  to  smaller length scales. This energy transference is physically realized through the phenomenon of vortex stretching where the eddies at the length scale $\ell_1< \ell_{0}$ break down into smaller eddies at the length scale $\ell_2<\ell_1<\ell_0$ and such  \emph{cascade} of energy transference proceeds until we reach the Kolmogorov dissipation scale  $\ell_D<\cdots <\ell_2<\ell_1<\ell_0$. Below this length scale $\ell_D$, the kinetic energy coming from larger scales is  ultimately dissipated by the direct action of the molecular viscosity and thus we will call the \emph{inertial range} the interval of length scales $]\ell_D,\ell_0[$ where the kinetic energy is transferred.\\

%%%%%%%%%%%%%%%%%%%%%%%%%%%%%%%%%
In order to state the \emph{Kolmogorov dissipation law} we need to introduce some terminology: let $\varepsilon>0$ be the  energy dissipation rate which determines the amount of energy lost by the viscous forces (\emph{i.e.} below the Kolmogorov dissipation scale $\ell_D$)  in the turbulent flow and which is given as an average of the gradient of the velocity (see (\ref{varepsilon,U}) for a precise definition). Define $U>0$ as the fluid characteristic velocity which is given by an average of the velocity  (see (\ref{varepsilon,U}) below) and  consider $L\geq \ell_0$ to be the fluid characteristic length which is related to the domain where we study the Kolmogorov dissipation law. Now, from the quantities $U$ and $L$ and using the viscosity parameter $\nu$, we can define the Reynolds number $Re$ by the expression (see \cite{Const,DoerFoias}):
\begin{equation}\label{Reynolds}
Re=\frac{UL}{\nu}.
\end{equation}  
At this point it is worth to do the following comments. Remark first that since the $U$ has a physical dimension of $length/time$,  $\nu$ has a physical dimension of $length^2/time$ and $L$ has a physical dimension of $length$,  the Reynolds number $Re$ has not any physical dimension. The second remark concerns the fact that the definition of the Reynolds number is \emph{not universal}  and  in the K41 theory this number is commonly defined though the energy input scale $\ell_0$ instead of the fluid characteristic length $L \geq \ell_0$ (see \emph{e.g.}  \cite{DoerFoias,OttoRamos}). But, in our definition of the Reynolds number   (\ref{Reynolds}),  we shall consider the length $L$ and this choice is motivated by the fact that in our model the fluid characteristic length $L$, defined in formula (\ref{LC}) below and which depends only on the data of our model, can be set large enough and thus we obtain large values for the Reynolds number $Re$.  \\

Then the Kolmogorov dissipation law states that in a turbulent setting \emph{i.e.} when $Re$ is large enough, the energy dissipation rate $\varepsilon$ can be estimated from above and from below in the following manner
\begin{equation}\label{dissipation_law}
c_1\frac{U^3}{\ell_0}\leq \varepsilon \leq c_2 \frac{U^3}{\ell_0}, \quad \text{when}\quad Re\gg1, 
\end{equation} 
where, $c_1,c_2>0$ are constants that \emph{do not depend} of the Reynolds number $Re$ (see \cite{Kolm1,Kolm2,Obuk}). Note that estimate (\ref{dissipation_law}) although often observed in many experiments (see \cite{Houel,Tennekes,Wilcox}) has a purely phenomenological explanation and its mathematical study constitutes a major challenge.\\

The general aim of this article is to study this dissipation law in the setting of the damped Navier-Stokes equations (\ref{damped_alpha_kappa_N-S}) defined over the whole space $\Rt$.\\

Several remarks are necessary before we proceed to the statement of the theorems. First let us stress that the ``characteristic length'' $L$ is intended to be the largest length scale where we shall observe a turbulent behavior: if we consider a fluid in a bounded domain $\Omega \subset \Rt$, this length scale is naturally linked to the size of the domain $\Omega$ and if we consider a periodic fluid on the cube $[0,L]^3$ it is reasonable to set the characteristic length as the period $L$. However since we work here with a fluid defined in the whole space $\Rt$ a natural approach to the characteristic length is not a completely trivial question. In this article we will define this length $L$ as a function of the external force (see formula (\ref{LC}) below) motivated by the fact that if any turbulence is observed, then it should appear where the action of the external force is actually relevant. \\

Our second remark is related to the force $\fe\in L^2(\Rt)\cap \dot{H}^{-1}(\Rt)$ which is divergence free and time-independent. This stationarity assumption is a simplification of the model since if we consider a time-dependent force then we will need to find a sufficiently large time interval in which the fluid is turbulent and this is a highly non-trivial issue. To solve this problem and since we have a constant supply of energy due to the external force we will study the turbulent behavior of the fluid by considering large time averages.\\ 

The third remark concerns the energy dissipation rate $\varepsilon$ and the fluid characteristic velocity $U$. In the articles \cite{DoerFoias,FMRT1,OttoRamos}, where it is considered a periodic fluid on the cube $[0,L]^3$, it is suggested to define these quantities by the following long-time averages: $\varepsilon =\nu \ds{\limsup_{T\longrightarrow +\infty}} \frac{1}{T}\int_{0}^{T}\Vert \vu(t,\cdot)\Vert^{2}_{\dot{H}^1}\frac{dt}{L^3}$ and $U=\left( \ds{\limsup_{T\longrightarrow +\infty}}  \frac{1}{T}\int_{0}^{T}\Vert \vu(t,\cdot)\Vert^{2}_{L^2}\frac{dt}{L^3}\right)^{\frac{1}{2}}$, 
where $L>0$ is the period. But, as we consider here a fluid in the whole space ( where a natural definition of the fluid characteristic length is not trivial)  in order to define the average quantities $\varepsilon, U$ we propose to consider the energy input scale $\ell_0>0$ which is a fix datum of our model, and thus we shall define:
\begin{equation}\label{varepsilon,U}
\varepsilon =\nu\limsup_{T\longrightarrow +\infty} \frac{1}{T}\int_{0}^{T}\Vert \vu(t,\cdot)\Vert^{2}_{\dot{H}^1}\frac{dt}{\ell^3_{0}} \quad\text{and}\quad U=\left(\limsup_{T\longrightarrow +\infty}  \frac{1}{T}\int_{0}^{T}\Vert \vu(t,\cdot)\Vert^{2}_{L^2}\frac{dt}{\ell^3_{0}}\right)^{\frac{1}{2}},
\end{equation}  Our first task is to give a sense to these quantities in the general framework of Leray's weak solutions\footnote{The use of the space $(L^{2}_{t})_{loc}$ in the time variable is related to the presence of the stationary external force.} $\vu\in L^{\infty}_{t}L^{2}_{x}\cap (L^{2}_{t})_{loc}\dot{H}^{1}_{x}$. The energy inequality verified by $\vu$ allows us to prove that $\varepsilon <+\infty$ (see the Appendix \ref{AppendixA} for a short proof of this fact) but the fact that $U<+\infty$ is highly non trivial in $\mathbb{R}^{3}$.  Observe that in the periodic setting $[0,L]^{3}$ we have at our disposal the following Poincaré inequality $\ds{\Vert \vu(t,\cdot) \Vert_{L^2} \leq \frac{ L}{2\pi} \Vert \vec{\nabla} \otimes \vu(t,\cdot)\Vert_{L^2}}$ which combined with the energy inequality satisfied by $\vu$ allows us to deduce that $U<+\infty$ (see  the Appendix \ref{AppendixB}). However, if we want to consider the quantity  $U$ over the whole space $\Rt$ we face with some important technical problems as we can not use the Poincaré inequality anymore and for a generic Leray's weak solution of the classical Navier-Stokes equations we can not assure (to the best of our knowledge) that the characteristic velocity $U$ is a finite quantity. Indeed,  from the classical energy inequality
\begin{equation}\label{energy-classical}
 \Vert \vu(t,\cdot)\Vert^{2}_{L^2}+2\nu\int_{0}^{t} \Vert 
  \vu(s,\cdot)\Vert^{2}_{\dot{H}^{1}}ds\leq \Vert \vu_0 \Vert^{2}_{L^2}+2\int_{0}^{t} \langle \fe, \vu(s,\cdot)  \rangle_{\dot{H}^{-1}\times \dot{H}^1}  ds,
\end{equation} we may prove the following control in time (see the details in the Appendix \ref{AppendixC}):  
\begin{equation}\label{control1}
\Vert \vu(t,\cdot)\Vert^{2}_{L^2}\leq \Vert \vu_0\Vert^{2}_{L^2}+\frac{t}{\nu}\Vert \fe \Vert^{2}_{\dot{H}^{-1}},
\end{equation}
but, when we apply a long-time average in this inequality  we get 
$$U^2=\limsup_{T\longrightarrow +\infty}\frac{1}{T}\int_{0}^{T}\Vert \vu(t,\cdot)\Vert^{2}_{L^2}\frac{dt}{\ell^{3}_{0}}\leq \limsup_{T\longrightarrow +\infty} T \frac{\Vert \fe \Vert^{2}_{\dot{H}^{-1}}}{2\nu \ell^{3}_{0}}=+\infty,$$ 
and we do not know if it is possible to obtain a better control in time for the quantity $\Vert \vu(t,\cdot)\Vert^{2}_{L^2}$ in order to assure that we actually  have $U<+\infty$. As a consequence of this fact, in the setting of the classical Navier-Stokes equations over the whole space $\Rt$, the study of the Kolmogorov dissipation law (\ref{dissipation_law}) could be potentially ill-posed and to overcome this problem we propose in this article to introduce a damping term which will allow us to prove that for any Leray's weak solution of equations (\ref{damped_alpha_kappa_N-S}) we have $U<+\infty$.\\

The fourth remark is related to the Reynolds number $Re$: note that the turbulent regime of a fluid can be characterized by the condition $Re\gg 1$ but this is an a posteriori approach as from formula (\ref{Reynolds}) the knowledge of the characteristic velocity $U$ is needed in order to determine the number $Re$. Since we have to handle carefully the definition of $U$, we shall use instead the \emph{Grashof number} $Gr$ (see (\ref{Grashof}) for a precise definition) which do not depend on $U$ and are equivalent to the Reynolds number.\\

The next remark  is about inequalities (\ref{dissipation_law}). Indeed, even if in the periodic framework we can assure that $\varepsilon, U<+\infty$, only the upper estimate $\ds{\varepsilon \leq c_2 \frac{U^3}{\ell_0}}$ is known under some technical conditions and the lower bound $\ds{c_1\frac{U^3}{\ell_0}\leq \varepsilon}$ is still an open problem (see the articles \cite{Const,DasGru,DoerFoias,FMRT1,FMRT2,Vigneron}). In the case of the whole space $\Rt$, the damping term  and the \emph{particular} characteristic length $L$ (see (\ref{LC}) below) will play an interesting role as they will help us to prove the lower and the upper bounds: $\ds{c_1 \frac{U^3}{L} \leq  \varepsilon \leq c_2 \frac{U^3}{L}}$, which are inequalities of the same type as in the Kolmogorov dissipation law  (\ref{dissipation_law}).\\

Our last remark focuses on the damping term $-\alpha P_\kappa(\vu)$ which is given by a frequency truncation as explained in expression (\ref{op_P}) below. In the deterministic study of turbulence it is quite natural to work in the Fourier variable (see the lecture notes \cite{Houel}, \cite{JacTab} and the Ph.D. thesis \cite{OscarTesis}) and thus the introduction of a cut-off function in the frequency level seems to be a well suited damping term. However we will show that for a particular set of the parameters $\alpha, \kappa >0$, the presence of the term $-\alpha P_\kappa(\vu)$ annihilates the turbulence. In order to highlight this phenomenon we will use a third characterization of the turbulence given by the \emph{Taylor length scale} $\ell_{T}$ which is defined in formula (\ref{relation_taylor_injection}) below. Indeed, in a turbulent regime we should have $Re\gg 1$, $Gr\gg 1$ and $\ell_{T}\ll \ell_{0}$ (see \cite{DoerFoias}, \cite{JacTab}), but due to the action of the damping term $-\alpha P_\kappa(\vu)$, and even if the Reynolds and Grashof numbers are large,  we will show that we actually have the equivalence $\ell_{T}\simeq \ell_{0}$ from which we can deduce that the fluid is \emph{not} in a turbulence regime.\\ 

As a conclusion, we can see that although the additional damping term is essential to assure that $U<+\infty$ and for proving the inequalities (\ref{dissipation_law}), the frequency truncation (which is in some sense a natural approach) introduced by the operator $-\alpha P_\kappa(\vu)$ is probably too strong and reduces significatively the effect of the expected turbulence. Let us finish observing that other damping terms can be considered in order to obtain $U<+\infty$ (see for example the article \cite{OJarrin}) but the complete study of the Kolmogorov dissipation law remains a challenging open problem.\\

The plan of the article is the following: in Section \ref{Secc_StatementResults} we will first introduce some notation and then we will state the theorems. In Section \ref{Sec_damped_N-S}, we will prove the existence of Leray's weak solutions for the damped Navier-Stokes equations (\ref{damped_alpha_kappa_N-S}) and in Section \ref{Secc_Ubounded} we will see how to deduce that the characteristic velocity $U$ is a finite quantity. In Section \ref{Secc_GrashofReynolds} we will prove that the Grashof numbers are equivalent to the Reynolds number and in Section \ref{Sec_Kolmogorov}, the Kolmogorov dissipation law (\ref{dissipation_law}) will be obtained for the damped equations (\ref{damped_alpha_kappa_N-S}). Finally, in Section \ref{Sec_non_turbulent}, we will show, by proving a general theorem and giving a precise example, how the extra damping term disturbes the effect of the turbulence.
%%%%%%%%%%%%%%%%%%%%%%%%%%%%%%%%%%%%%%
\section{Statement of the results}\label{Secc_StatementResults}
We give now the definition of the damping term $-\alpha P_\kappa(\vu)$ introduced in equations (\ref{damped_alpha_kappa_N-S}): for $\alpha,\kappa>0$ two positive parameters, this operator is defined in the Fourier variable by
\begin{equation}\label{op_P}
\widehat{-\alpha P_\kappa(\vu)}(t,\xi)=-\alpha\mathds{1}_{\vert \xi \vert<\kappa}(\xi)\widehat{\vu}(t,\xi),
\end{equation} 
where $\alpha>0$ is a damping parameter, $\kappa>0$ is a cut-off frequency and $\widehat{\vu}$ denotes the Fourier transform in the space variable of the function $\vu$.  For the time being we do not impose any restriction on the parameters $(\alpha, \kappa)$, but we will see later on that some conditions are required.\\

Our first theorem studies the existence of Leray's weak solutions of the damped Navier-Stokes equations (\ref{damped_alpha_kappa_N-S}).
%%%%%%%%%%%%%%%%%%%%%%%%%%%%%%%%%%%%%%
\begin{Theoreme}\label{Existence} 
Let $\vu_0\in L^2(\Rt)$ be a divergence-free initial data and consider $\fe\in L^2(\Rt)\cap\dot{H}^{-1}(\Rt)$ be a divergence-free, stationary external force. Then, for all $\alpha,\kappa>0$ there exists a function $\vu_{(\alpha,\kappa)} \in L^{\infty}(]0,+\infty[,L^2(\Rt))\cap L^{2}_{loc}(]0,+\infty[,\dot{H}^1(\Rt))$ which is a weak solution of the damped Navier-Stokes equations (\ref{damped_alpha_kappa_N-S}). Moreover, this solution satisfies the following energy inequality: 
\begin{equation}\label{energy_ineq_alpha_model}
\begin{split}
\Vert \vu_{(\alpha,\kappa)}(t,\cdot)\Vert^{2}_{L^2}+2\nu\int_{0}^{t} \Vert \vu_{(\alpha,\kappa)}(s,\cdot)\Vert^{2}_{\dot{H}^1}ds\leq \Vert \vu_0 \Vert^{2}_{L^2}+2\int_{0}^{t} \langle \fe, \vu_{(\alpha,\kappa)}(s,\cdot)\rangle_{\dot{H}^{-1}\times \dot{H}^1} ds \\
-2\alpha \int_{0}^{t}\Vert P_\kappa(\vu_{(\alpha,\kappa)})(s,\cdot) \Vert^{2}_{L^2}ds.
\end{split}
\end{equation}
\end{Theoreme}
%%%%%%%%%%%%%%%%%%%%%%%%%%%%%%%%%%%%%%
From now on, we will simply denote by $\vu$ the weak solution $\vu_{(\alpha,\kappa)}$ obtained in the previous theorem.\\

As pointed out in the introduction, the presence of the additional damping term allow us to prove that the characteristic velocity $U$ is actually a bounded quantity:
%%%%%%%%%%%%%%%%%%%%%%%%%%%%%%%%%%%%%%
\begin{Theoreme}\label{Th-control-in-time} Within the framework of Theorem \ref{Existence}, for all $\alpha,\kappa>0$, weak solutions $\vu \in L^{\infty}_{t}L^{2}_{x}\cap (L^{2}_{t})_{loc}\dot{H}^{1}_{x}$ of the damped Navier-Stokes equations  (\ref{damped_alpha_kappa_N-S}) satisfies the following estimate: for $\beta=\min\left( 2\alpha, \nu\kappa^2\right)$ and for all time $t\in [0,+\infty[$ we have,
\begin{equation}\label{control_in_time_Th_2}
\Vert \vu(t,\cdot)\Vert^{2}_{L^2}\leq e^{-\beta t}\Vert \vu_0 \Vert^{2}_{L^2}+c\frac{\Vert \fe \Vert^{2}_{\dot{H}^{-1}}}{\nu \beta} \left( 1-e^{-\beta t}\right),
\end{equation}
from which we can deduce the control
\begin{equation}\label{estim-U}	
U^2= \limsup_{T \longrightarrow +\infty} \frac{1}{T} \int_{0}^{T} \Vert \vu(t,\cdot)\Vert_{L^{2}}^{2} \frac{dt}{\ell^{3}_{0}} \leq c \frac{\Vert \fe \Vert^{2}_{\dot{H}^{-1}}}{\nu \beta \ell^{3}_{0}} <+\infty. 
\end{equation}

\end{Theoreme}
%%%%%%%%%%%%%%%%%%%%%%%%%%%%%%%%%%%%%%
Estimate (\ref{control_in_time_Th_2}) is of course the key to obtain that $U<+\infty$. Observe in particular that if $\alpha\to 0$ or if $\kappa\to 0$ then, due to the presence of the term $\frac{\Vert \fe \Vert^{2}_{\dot{H}^{-1}}}{\nu \beta}$ and from the definition of $\beta$, we can not longer deduce that $U<+\infty$. \\

The two previous results were quite general and in order to continue we need to fix the parameters of our model. First, we will consider a stationary and divergence-free external force $\fe\in L^2(\Rt)\cap \dot{H}^{-1}(\Rt)$ such that
\begin{equation}\label{loc_fe}
supp\,\left( \widehat{\fe}\right) \subset \left\lbrace \xi \in \Rt: \frac{1}{10 \theta \ell_0} \leq \vert \xi \vert \leq \frac{1}{\theta\ell_0}\right\rbrace,
\end{equation} 
where $\ell_0>0$ is an energy input scale that will be fixed from now on and $ \theta \geq 1$ is a technical parameter which does not depend on any physical parameter of our model and which will be useful later on. Remark that  this Fourier-support condition is often considered in the litterature to represent  the fact that the kinetic energy is only introduced in the fluid by the external force $\fe$ at the length scale of the order  $\ell_0$ and thus at the frequency of order $\frac{1}{\ell_0}$ (see also the articles \cite{Chesk}, \cite{Const2}, \cite{Const}, \cite{FMRT1}, \cite{FMRT2} and \cite{OttoRamos} for similar conditions).\\

We continue with the definition of the characteristic length $L$ and following an idea suggested by the lecture notes \cite{Const}, we shall define this quantity by the expression
\begin{equation}\label{LC}
L=\frac{\ell_0}{\gamma},
\end{equation} 
where $\ell_0$ is the energy input scale given by the external force $\fe$ and $\gamma>0$ is a parameter given as follows: since $\widehat{\fe}$ is localized at the frequencies $\frac{ 1}{10\theta  \ell_0}\leq \vert \xi \vert\leq \frac{1}{\theta\ell_0}$, by the Bernstein inequalities  there exists  a constant $c_0>0$ such that we have  the inequality (recall that $\theta \geq 1$)
\begin{equation}\label{c0}
\Vert \fe \Vert_{L^{\infty}}\leq c_0\left(\frac{1}{\theta\ell_0}\right)^{\frac{3}{2}}\Vert \fe \Vert_{L^2}\leq c_0\left(\frac{1}{\ell_0}\right)^{\frac{3}{2}}\Vert \fe \Vert_{L^2},
\end{equation} 
and from this inequality we define the parameter $\gamma>0$ by 
\begin{equation}\label{gamma}
\gamma =\frac{\Vert \fe \Vert_{L^{\infty}}}{c_0 \left(\frac{1}{\ell_0}\right)^{\frac{3}{2}}\Vert \fe \Vert_{L^2}}\leq 1. 
\end{equation} 
Remark that by formulas (\ref{LC}) and (\ref{gamma}) we have that the characteristic length $L$  depends only of the external force $\fe$. \\

Now, as said in the introduction, we will characterize the turbulent regime by a condition on the Grashof number $Gr$ which are given by
\begin{equation}\label{Grashof}
Gr=\frac{FL^3}{\nu^2},
\end{equation}
where 
\begin{equation}\label{F}
\ds{F= \frac{\Vert \fe \Vert_{L^2}}{\ell^{\frac{3}{2}}_{0}}},
\end{equation} 
is an average in spatial variable  of function $\fe$. Note that the Grashof number $Gr$ depend essentially on the external force $\fe$ and not on the characteristic velocity $U$. See the article \cite{DoerFoias} and the book \cite{Tennekes} for more details about the Grashof number.\\

We will prove that the Grashof number $Gr$ is \emph{equivalent} to the Reynolds number $Re$ and in order to study  this fact  we need to introduce the following parameter: consider the number $G_0>0$ defined as: 
\begin{equation}\label{Grashof_zeo}
G_0=\frac{\Vert \fe \Vert_{L^{\infty}}\ell^{3}_{0}}{\nu^2},
\end{equation}
and we have  the following relationship
\begin{equation}\label{estimate_Grashof_G0}
Gr=\frac{G_0}{c_0 \gamma^4},
\end{equation}
where the parameter $\gamma>0$ is given in formula (\ref{gamma}) and the numerical  constant $c_0>0$ is given in expression (\ref{c0}). In this identity we may see that we can obtain large values of the Grashof number $Gr$ by fixing the number $G_0$ by letting the parameter $\gamma$ to be small enough.\\ 

We fix now the parameters of the damping operator $-\alpha P_{\kappa}$ given in (\ref{op_P}):  for $\alpha$ and for the cut-off frequency $\kappa>0$ we will consider
\begin{equation}\label{alphakappa}
\alpha= \frac{\nu}{\ell^{2}_{0}},\quad \mbox{and}\quad\kappa= \frac{1}{20\, \theta\ell_0}.
\end{equation} 
We have then the following result:
%%%%%%%%% %%%%%%%%% %%%%%%%%% %%%%%%%%%%
\begin{Theoreme}\label{Proposition} 
Let $\fe\in L^2(\Rt)\cap \dot{H}^{-1}(\Rt)$ be a stationary and divergence-free external force which satisfies the frequency localization (\ref{loc_fe}) for some $\ell_0>0$ and some $\theta \geq 1$ large enough  and consider $\vu$ a solution of equation (\ref{damped_alpha_kappa_N-S}) with a damping term $-\alpha P_{\kappa}$ such that $\alpha,\kappa$ satisfy (\ref{alphakappa}). Let $Re$ and $Gr$ be the Reynolds and the Grashof number given by expressions (\ref{Reynolds}) and (\ref{Grashof}) and let $G_0>0$ be the number given in (\ref{Grashof_zeo}). Then we have  
\begin{equation}\label{estimate_reynolds_grashof}
\mathfrak{a}_{1,G_{0}}\, Re^{2}  \leq  Gr \leq  \mathfrak{a}_{2,G_{0}}\, Re^{2},
\end{equation}
where the constants $0<\mathfrak{a}_{1,G_{0}}\leq \mathfrak{a}_{2,G_{0}}$  depend on $G_{0}$.
\end{Theoreme}
\begin{Remarque}
	\emph{The definition of the parameter $\alpha$ above is essentially meant to obtain dimensionless constants,  whereas the upper bound of the condition $1\leq \theta$ is related to the frequency cut-off $\kappa$ as we need to separate the action of the force from the action of the damping term.}
\end{Remarque}

We are now ready to state one of the main results of this article in which we will prove the Kolmogorov dissipation law (\ref{dissipation_law}) for the damped Navier-Stokes equations (\ref{damped_alpha_kappa_N-S}):
%%%%%%%%% %%%%%%%%% %%%%%%%%% %%%%%%%%%%
\begin{Theoreme}\label{Main_Result} Let $\nu>0$ be the fluid viscosity parameter and let  $\fe\in L^2(\Rt)\cap \dot{H}^{-1}(\Rt)$ be a stationary and divergence-free external force which satisfies the frequency localization (\ref{loc_fe}) with $\ell_0>0$ a fixed energy input scale. 
	
	Then,  if the Grashof number $Gr$ given in (\ref{Grashof}) satisfies the condition
$\ds{ Gr \geq \frac{4  \mathfrak{a}_{2,G_0} }{c^{2}_{0}} \frac{G^{2}_{0}}{\gamma^4}}$,  where the constant $\mathfrak{a}_{2,G_0}$ is given by (\ref{estimate_reynolds_grashof}) above, there exists two constants $0<\mathfrak{b}_{1,G_0}\leq \mathfrak{b}_{2,G_0}$, which depend on $G_0>0$, such that we have the Kolmogorov dissipation law
$$\mathfrak{b}_{1,G_0}\frac{U^3}{L} \leq \varepsilon \leq \mathfrak{b}_{2,G_0} \frac{U^3}{L},$$
where $L>0$ is the fluid characteristic length defined in (\ref{LC}) and the quantities $\varepsilon$ and $U$ given in (\ref{varepsilon,U}) are built from a weak solution $\vu \in L^{\infty}([0,+\infty[,L^{2}(\Rt))\cap L^{2}_{loc}([0,+\infty[,\dot{H}^1(\Rt))$ obtained in Theorem \ref{Existence}.
\end{Theoreme}
\begin{Remarque}
	\emph{Recall that by (\ref{estimate_Grashof_G0}) we have the identity $\ds{Gr= \frac{G_0}{c_0 \gamma^4}}$, and then the condition on the Grashof $Gr$ above is satisfied when the number $G_0$ is small enough and verifies: $\ds{\frac{4 \mathfrak{a}_{2,G_0} }{c_0}} G_0 \leq 1$.}
\end{Remarque}

%%%%%%%%% %%%%%%%%% %%%%%%%%% %%%%%%%%%%

Our last theorem studies more carefully the framework introduced until now and we will see that, although we were able to obtain the Kolmogorov dissipation law for turbulent fluids with Theorem \ref{Main_Result}, all this framework is actually non turbulent. To be more precise, we need to introduce the \emph{Taylor scale} $\ell_T$ which is defined as follows:  
\begin{equation}\label{Taylor_scale}
\ell_T=\left( \frac{\nu U^2}{\varepsilon}\right)^{\frac{1}{2}},
\end{equation}
where $\nu$ is the viscosity parameter and $U$ and $\varepsilon$ are the characteristic velocity and the energy dissipation rate, respectively. See the article \cite{DoerFoias} and the book \cite{JacTab} for more references about the Taylor scale. This length scale, also called the \emph{turbulence length scale}, is commonly used to characterize the turbulent regime:  according to the Kolmogorov dissipation law for large Reynolds number  $Re\gg1$ we should have the equivalence:
\begin{equation}\label{relation_taylor_injection}
\ell_T \simeq \frac{1}{\sqrt{Re}}\ell_0,
\end{equation}
which can be expressed (due to formula (\ref{estimate_reynolds_grashof})) in terms of the  Grashof number  as $\ds{\ell_T \simeq \frac{1}{(Gr)^{\frac{1}{4}}}\ell_0}$. Thus, in a turbulent regime (where $Re\gg 1$ or $Gr\gg 1$) we should have
\begin{equation}\label{Taylor-turbulence}
\ell_T \ll\ell_0,
\end{equation}
and the study of this relationship is exactly the purpose of the following theorem. 
\begin{Theoreme}\label{Th-non-turb} Under the same hypotheses of Theorem \ref{Main_Result}, consider $\ell_T>0$ the Taylor scale defined in equation (\ref{Taylor_scale}). Then, there exists two constants $0<\mathfrak{c}_{1,G_0}\leq \mathfrak{c}_{2,G_0}$, which depend on $G_0>0$, such that we have the estimates:
\begin{equation}\label{estimate_taylor_2}
\mathfrak{c}_{1,G_0}\, \ell_0 \leq \ell_T \leq \mathfrak{c}_{2,G_0}\, \ell_0. 
\end{equation}
\end{Theoreme}	

We can thus see that, even if the Reynolds number is large and even if the Kolmogorov dissipation law is satisfied, the particular choice of the parameters $\alpha$ and $\kappa$ used in the previous theorems can annihilate the turbulence. Maybe other values of these parameters can be considered, but if there is any interference at the Fourier level between the damping term considered here and the external force this must be treated very carefully and it is, to the best of our knowledge, an interesting open problem. 

%%%%%%%%%%%%%%%%%%%%%%%%%%%%%%%%%%%%%%%%%%%%%%%%%%%%%%%%%%%%%%%%%%%%%%%%%%%%%%%%%%%%%%%%%%%%%%%%%%%%%%%%%%%%%%%%%%%%%%%%%%%%%%%%%%%%%%%%%%%%%%%%%%%%%%%%%%
\section{Theorem \ref{Existence}: existence for the damped Navier-Stokes equations}\label{Sec_damped_N-S}
The proof of the existence of weak solutions for equation (\ref{damped_alpha_kappa_N-S}) is rather straightforward as it follows essentially the same lines than the classical framework. Let $\phi \in \mathcal{C}^{\infty}_{0}(\Rt)$ be a positive function such that $\ds{\int_{\Rt}\phi(x)dx=1}$, for $\delta>0$ and for $\phi_\delta(x)=\frac{1}{\delta^3}\phi \left( \frac{x}{\delta}\right)$, we will solve the following integral equation 
\begin{eqnarray}\label{NS_regularised_point_fix}\nonumber
\vu(t,x)&=&h_{\nu t}\ast \vu_0(x) +\int_{0}^{t}h_{\nu(t-s)}\ast\fe(x)ds-\int_{0}^{t}h_{\nu(t-s)}\ast(\P(([\phi_\delta\ast\vu]\cdot \vec{\nabla})\vu)(s,x)ds\\ 
& &-\alpha\int_{0}^{t}h_{\nu(t-s)}\ast P_\kappa(\vu)(s,x)ds,
\end{eqnarray}
in the space  $L^{\infty}([0,T],L^2(\Rt))\cap L^2([0,T],\dot{H}^1(\Rt))$ provided with the norm $\Vert \cdot \Vert_T=\Vert \cdot \Vert_{L^{\infty}_{t}L^{2}_{x}}+\sqrt{\nu}\Vert \cdot \Vert_{L^{2}_{t}\dot{H}^{1}_{x}}$. We have then:
\begin{eqnarray*}
\Vert \vu \Vert_T &\leq &\underbrace{\left\Vert h_{\nu t}\ast \vu_0 +\int_{0}^{t}h_{\nu(t-s)}\ast\fe(\cdot)ds\right\Vert_T}_{(1)}  + \underbrace{\left\Vert \int_{0}^{t}h_{\nu(t-s)}\ast(\P(([\phi_\delta\ast\vu]\cdot \vec{\nabla})\vu)(s,\cdot)ds\right\Vert_T}_{(2)} \\
& & +\underbrace{\alpha  \left\Vert \int_{0}^{t}h_{\nu(t-s)}\ast P_\kappa(\vu)(s,\cdot)ds \right\Vert_T}_{(3)}.
\end{eqnarray*}
Terms $(1)$ and $(2)$ are classical to estimate. Indeed, by \cite{PGLR1}, Theorem $12.2$, page $352$,  we have $\Vert h_{\nu t}\ast \vu_0\Vert_{T}\leq c \Vert \vu_0 \Vert_{L^2}$ and  
\begin{equation}\label{estimate_force_th_existence}
\left\Vert \int_{0}^{t}h_{\nu(t-s)}\ast\fe(\cdot)ds \right\Vert_T \leq C(\frac{1}{\sqrt{\nu}}+\sqrt{T}) \Vert \fe \Vert_{L^{2}_{t}H^{-1}_{x}},
\end{equation} 
but since $\fe \in L^{2}(\Rt)\cap \dot{H}^{-1}(\Rt)$ we have: $\Vert \fe \Vert_{L^{2}_{t}H^{-1}_{x}}\leq \ds{\Vert \fe \Vert_{L^{2}_{t} \dot{H}^{-1}_{x} } \leq \sqrt{T} \Vert \fe \Vert_{L^{\infty}_{t} \dot{H}^{-1}_{x}} \leq \sqrt{T} \Vert \fe \Vert_{\dot{H}^{-1}}}$, and thus, for term $(1)$ above we can write  
\begin{equation}\label{(1)}
\left\Vert h_{\nu t}\ast \vu_0 +\int_{0}^{t}h_{\nu(t-s)}\ast\fe(\cdot)ds\right\Vert_T  \leq  c\Vert \vu_0 \Vert_{L^2}+ C(\frac{1}{\sqrt{\nu}}+\sqrt{T}) \sqrt{T}\Vert \fe \Vert_{\dot{H}^{-1}}.
\end{equation} 
For term $(2)$ we have (see \cite{PGLR1}, Theorem $12.2$ for the details):
\begin{equation}\label{(2)}
\left\Vert \int_{0}^{t}h_{\nu(t-s)}\ast(\P(([\phi_\delta\ast\vu]\cdot \vec{\nabla})\vu)(s,\cdot)ds\right\Vert_T \leq C (\frac{1}{\sqrt{\nu}}+1) \sqrt{T}\delta^{-\frac{3}{2}} \Vert \vu \Vert_T\,\Vert \vu \Vert_T.
\end{equation}
We only need now to study the quantity $(3)$ and since we have $\widehat{P_\kappa(\vu)}(t,\xi)=\mathds{1}_{\vert \xi \vert <\kappa}(\xi)\widehat{\vu}(t,x)$ then by the Plancherel identity we can write  
\begin{eqnarray*}
\Vert P_\kappa(\vu) \Vert_{L^{2}_{t}H^{-1}_{x}}&=& \left( \int_{0}^{T}\left\Vert (1+\vert \xi \vert^2)^{-\frac{1}{2}}\widehat{P_\kappa(\vu)}(t,\cdot)\right\Vert^{2}_{L^2}dt\right)^{\frac{1}{2}}= \left( \int_{0}^{T}\left\Vert (1+\vert \xi \vert^2)^{-\frac{1}{2}}\left(\mathds{1}_{\vert \xi \vert<\kappa}(\xi) \widehat{\vu}(t,\cdot)\right)\right\Vert^{2}_{L^2}dt\right)^{\frac{1}{2}}\\
&\leq & \left( \int_{0}^{T}\left\Vert \widehat{\vu}(t,\cdot)\right\Vert^{2}_{L^2}dt\right)^{\frac{1}{2}}\leq \sqrt{T}\Vert \widehat{\vu}\Vert_{L^{\infty}_{t}L^{2}_{x}}=\sqrt{T}\Vert\vu\Vert_{L^{\infty}_{t}L^{2}_{x}}\leq \sqrt{T}\Vert \vu \Vert_T.
\end{eqnarray*}
Now, substituting $\fe$ by $\alpha P_\kappa(\vu)$ in inequality (\ref{estimate_force_th_existence}) above we have 
$$ \alpha  \left\Vert \int_{0}^{t}h_{\nu(t-s)}\ast P_\kappa(\vu)(s,\cdot)ds \right\Vert_T \leq \alpha C \left( \frac{1}{\sqrt{\nu}}+\sqrt{T}\right)\Vert P_\kappa(\vu) \Vert_{L^{2}_{t}H^{-1}_{x}} \leq \alpha C \left( \frac{1}{\sqrt{\nu}}+\sqrt{T}\right)\sqrt{T}\Vert \vu \Vert_T,$$ and then we obtain the estimate 
\begin{equation}\label{(3)}
\alpha  \left\Vert \int_{0}^{t}h_{\nu(t-s)}\ast P_\kappa(\vu)(s,\cdot)ds \right\Vert_T \leq \alpha C\left( \frac{1}{\sqrt{\nu}}+\sqrt{T}\right)\sqrt{T}\Vert \vu\Vert_{T}.
\end{equation} 
Once we have inequalities (\ref{(1)}), (\ref{(2)}) and (\ref{(3)}),  for a  time $T>0$ small enough and for $\delta>0$, by the Banach contraction principle we obtain $\vu_\delta \in L^{\infty}([0,T],L^2(\Rt))\cap L^2([0,T],\dot{H}^1(\Rt))$ a local solution of  the equations (\ref{NS_regularised_point_fix}).\\

Now we will prove that this solution $\vu_\delta$ is global. Remark that the function $\vu_\delta$ satisfies the regularized equation 
\begin{equation}\label{NS-damped-regularized}
\partial_t\vu_\delta= \nu\Delta \vu_\delta-\P(([\phi_\delta\ast\vu_\delta]\cdot\vec{\nabla}) \vu_\delta)+\fe-\alpha P_\kappa(\vu_\delta),
\end{equation}
where all the terms belong to the space $L^2([0,T], \dot{H}^{-1}(\Rt))$ and then we can write 
\begin{eqnarray*}
\frac{d}{dt} \Vert \vu_\delta(t,\cdot)\Vert^{2}_{L^2} &=& 2\langle \partial_t \vu_\delta(t,\cdot),\vu_\delta(t,\cdot)\rangle_{\dot{H}^{-1}\times \dot{H}^1}\\
&=&-2\nu \Vert \vu_\delta(t,\cdot)\Vert^{2}_{\dot{H}^1}+2\langle \fe,\vu_\delta(t,\cdot)\rangle_{\dot{H}^{-1}\times \dot{H}^1}-2\alpha\langle P_\kappa(\vu_\delta)(t,\cdot), \vu_\delta(t,\cdot)\rangle_{\dot{H}^{-1}\times \dot{H}^1}.
\end{eqnarray*}
now since $\vu_{\delta}(t,\cdot)\in L^{2}(\Rt)$ then we have $P_\kappa(\vu_\delta)(t,\cdot)\in L^{2}(\Rt)$ and thus
\begin{eqnarray*}
\langle P_\kappa(\vu_\delta)(t,\cdot), \vu_\delta(t,\cdot)\rangle_{\dot{H}^{-1}\times \dot{H}^1} &=&\int_{\Rt}\widehat{P_\kappa(\vu_\delta)}(t,\xi) \cdot \widehat{\vu_\delta}(t,\xi) d\xi=\int_{\Rt} \left( \mathds{1}_{\vert \xi \vert<\kappa}(\xi)\widehat{\vu_\delta}(t,\xi)\right)\cdot \widehat{\vu_\delta}(t,\xi)d\xi\\ 
&=& \int_{\Rt} \left( \mathds{1}_{\vert \xi \vert<\kappa}(\xi)\widehat{\vu_\delta}(t,\xi)\right)\cdot \left(\mathds{1}_{\vert \xi \vert<\kappa}(\xi) \widehat{\vu_\delta}(t,\xi)\right) d\xi\\
&=&\Vert \widehat{P_\kappa(\vu_\delta)}(t,\cdot)\Vert^{2}_{L^2} = \Vert P_\kappa (\vu_\delta) (t,\cdot)\Vert^{2}_{L^2},
\end{eqnarray*}
and then we have
\begin{equation}\label{energ_eq_delta}
\frac{d}{dt} \Vert \vu_\delta(t,\cdot)\Vert^{2}_{L^2} = -2\nu \Vert \vu_\delta(t,\cdot)\Vert^{2}_{\dot{H}^1}+2\langle \fe,\vu_\delta(t,\cdot)\rangle_{\dot{H}^{-1}\times \dot{H}^1}-2\alpha \Vert P_\kappa(\vu_\delta)(t,\cdot)\Vert^{2}_{L^2}.
\end{equation} 
But $-2\alpha \Vert P_\kappa(\vu_\delta)(t,\cdot)\Vert^{2}_{L^2}$ is a negative quantity and we get
\begin{eqnarray*}
\frac{d}{dt} \Vert \vu_\delta(t,\cdot)\Vert^{2}_{L^2} &  \leq &-2\nu \Vert \vu_\delta(t,\cdot)\Vert^{2}_{\dot{H}^1} + 2\langle f,\vu_\delta(t,\cdot)\rangle_{\dot{H}^{-1}\times \dot{H}^1} \leq  -2\nu \Vert \vu_\delta (t,\cdot)\Vert^{2}_{\dot{H}^1} +\nu \Vert \vu_\delta(t,\cdot)\Vert^{2}_{\dot{H}^1}+\frac{1}{\nu}\Vert f \Vert^{2}_{\dot{H}^{-1}} \\
& \leq &  -\nu \Vert \vu_\delta (t,\cdot)\Vert^{2}_{\dot{H}^1}+\frac{1}{\nu}\Vert f \Vert^{2}_{\dot{H}^{-1}}.
\end{eqnarray*}
Finally, we integrate on the interval of time $[0,t]$  and  we obtain the following control
\begin{equation}\label{control_unif_delta}
\Vert \vu_\delta (t,\cdot)\Vert^{2}_{L^2}+\nu\int_{0}^{t} \Vert \vu_\delta(s,\cdot)\Vert^{2}_{\dot{H}^1}ds \leq \left(\Vert\vu_0\Vert^{2}_{L^2}+\frac{t}{\nu}\Vert f \Vert^{2}_{\dot{H}^{-1}}\right),
\end{equation} 
which allows us to extend the local solution $ \vu_\delta $ to the whole interval $[0,+\infty[$.\\
\\
We study now the convergence to a weak solution of equations (\ref{damped_alpha_kappa_N-S}). Indeed, by the  Rellich-Lions lemma (see \cite{PGLR1}, Theorem $12.1$) there exists a sequence of positive numbers $(\delta_n)_{n\in\mathbb{N}}$ and a function $\vu \in L^{2}_{loc}([0,+\infty[\times \Rt)$ such that the sequence $(\vu_{\delta_n})_{n\in \mathbb{N}}$  converges strongly to $\vu$ in $L^{2}_{loc}([0,+\infty[\times \Rt)$. Moreover,  this sequence converges  to $\vu$  in the weak$-*$ topology of the spaces  $L^{\infty}([0,T],L^2(\Rt))$ and $L^{2}([0,T],\dot{H}^1(\Rt))$ for all $T>0$. 
From these convergences we can deduce that the sequence $ \left(\P(([\phi_{\delta_n}\ast\vu_{\delta_n}]\cdot\vec{\nabla}) \vu_{\delta_n})\right)_{n\in \mathbb{N}}$ converges  to $\P((\vu\cdot\vec{\nabla}) \vu)$ in the weak$-*$ topology of the space  $(L^{2}_{t})_{loc}(H^{-\frac{3}{2}}_{x})$ and then,  in order to verify that the limit $\vu$ is a weak  solution of equation  (\ref{damped_alpha_kappa_N-S}) it remains to prove the convergence of the sequence  $\left( -\alpha P_\kappa(\vu_{\delta_n})\right)_{n\in \mathbb{N}}$ to $-\alpha P_\kappa(\vu)$. More precisely,  we will prove that the sequence  $(-\alpha P_\kappa(\vu_{\delta_n}))_{n\in\mathbb{N }}$ converges  to $-\alpha P_\kappa(\vu)$ in the weak$-*$ topology of the space $L^{2}([0,T],L^2(\Rt))$. Indeed, since we have  $\Vert \vu_{\delta_n}\Vert_{L^{2}([0,T],L^2(\Rt))}\leq \sqrt{T}\Vert \vu_{\delta_n}\Vert_{L^{\infty}([0,T],L^2(\Rt))}$ then by  inequality (\ref{control_unif_delta}) we get that the sequence $(\vu_{\delta_n})_{n\in\mathbb{N}}$ converges   to $\vu$ in the weak$-*$ topology of the space in $L^{2}([0,T],L^2(\Rt))$ and moreover, since $P_\kappa(\cdot)$ is strongly continuous in this space then we obtain the desired convergence. \\ 
\\
Finally, for the energy inequality we get back to identity (\ref{energ_eq_delta}) and we integrate each term of this identity:
$$\Vert \vu_{\delta_n}(t,\cdot)\Vert^{2}_{L^2}+2\nu \int_{0}^{t} \Vert \vu_{\delta_n}(s,\cdot)\Vert^{2}_{\dot{H}^1}ds = \Vert \vu_0\Vert^{2}_{L^2} +2 \int_{0}^{t}\langle \fe, \vu_{\delta_n}(s,\cdot)\rangle_{\dot{H}^{-1}\times \dot{H}^1}ds -2\alpha \int_{0}^{t}\Vert P_\kappa(\vu_{\delta_n})(s,\cdot)\Vert^{2}_{L^2}ds,$$
from which we obtain the inequality (\ref{energy_ineq_alpha_model}) by applying classical tools (see the book \cite{PGLR1}). Theorem \ref{Existence} is proven. \finpv 
%%%%%%%%%%%%%%%%%%%%%%%%%%%%%%%%%%%%%%
\section{Proof of Theorem \ref{Th-control-in-time}}\label{Secc_Ubounded}
Consider the functions $\vu_{\delta_n}$ which are solutions of the regularized equation (\ref{NS-damped-regularized}). Our starting point is then the equality (\ref{energ_eq_delta}):
$$ \frac{d}{dt} \Vert \vu_{\delta_n}(t,\cdot)\Vert^{2}_{L^2} =  -2\nu \Vert \vu_{\delta_n}(t,\cdot)\Vert^{2}_{\dot{H}^1}+2\langle \fe,\vu_{\delta_n}(t,\cdot)\rangle_{\dot{H}^{-1}\times \dot{H}^1}-2\alpha \Vert P_\kappa(\vu_{\delta_n})(t,\cdot)\Vert^{2}_{L^2}.$$ 
Since  $\fe\in \dot{H}^{-1}(\Rt)$ we have $ \ds{2\langle \fe,\vu_{\delta_n}(t,\cdot)\rangle_{\dot{H}^{-1}\times \dot{H}^1}\leq \frac{\Vert \fe \Vert^{2}_{\dot{H}^{-1}}}{\nu}+\nu \Vert \vu_{\delta_n}(t,\cdot)\Vert^{2}_{\dot{H}^1}}$  and then we get 
\begin{equation}\label{ineq_der_time_delta_2}
\frac{d}{dt} \Vert \vu_{\delta_n}(t,\cdot)\Vert^{2}_{L^2}\leq\frac{\Vert \fe\Vert^{2}_{\dot{H}^{-1}}}{\nu}  -\nu\Vert \vu_{\delta_n}(t,\cdot)\Vert^{2}_{\dot{H}^1} -2\alpha \Vert P_\kappa(\vu_{\delta_n})(t,\cdot)\Vert^{2}_{L^2}.
\end{equation}
Note that by the Plancherel identity we have  
\begin{eqnarray*}
-\nu\left\Vert \vu_{\delta_n}(t,\cdot)\right\Vert^{2}_{\dot{H}^1}& = & -\nu \left\Vert \vert \xi \vert \widehat{\vu_{\delta_n}}(t,\cdot)\right\Vert^{2}_{L^2}= -\nu \left\Vert \mathds{1}_{\vert \xi \vert<\kappa}  \vert \xi \vert \widehat{\vu_{\delta_n}}(t,\cdot)\right\Vert^{2}_{L^2}  -\nu \left\Vert \mathds{1}_{\vert \xi \vert\geq \kappa} \vert \xi \vert \widehat{\vu_{\delta_n}}(t,\cdot)\right\Vert^{2}_{L^2} \\
&\leq & -\nu \left\Vert \mathds{1}_{\vert \xi \vert\geq \kappa}  \vert \xi \vert \widehat{\vu_{\delta_n}}(t,\cdot)\right\Vert^{2}_{L^2}\leq -\nu \kappa^2\left\Vert \mathds{1}_{\vert \xi \vert\geq \kappa} \widehat{\vu_{\delta_n}}(t,\cdot)\right\Vert^{2}_{L^2},\\
\end{eqnarray*} 
moreover, always by the Plancherel identity, we obtain 
$$-\nu\left\Vert \vu_{\delta_n}(t,\cdot)\right\Vert^{2}_{\dot{H}^1} - 2\alpha \left\Vert P_\kappa(\vu_{\delta_n})(t,\cdot)\right\Vert^{2}_{L^2} \leq  -\nu \kappa^2 \left\Vert \mathds{1}_{\vert \xi \vert\geq \kappa}(\xi)  \widehat{\vu_{\delta_n}}(t,\cdot)\right\Vert^{2}_{L^2} -2\alpha\left\Vert \mathds{1}_{\vert \xi \vert <\kappa}(\xi)  \widehat{\vu_{\delta_n}}(t,\cdot)\right\Vert^{2}_{L^2}$$
\begin{eqnarray}\label{ineq_Plancherel}\nonumber
&\leq & -\min(2\alpha,\nu\kappa^2)\left( \left\Vert \mathds{1}_{\vert \xi \vert <\kappa}(\xi)  \widehat{\vu_{\delta_n}}(t,\cdot)\right\Vert^{2}_{L^2}+\left\Vert \mathds{1}_{\vert \xi \vert \geq \kappa}(\xi)  \widehat{\vu_{\delta_n}}(t,\cdot)\right\Vert^{2}_{L^2}\right)\\ 
&\leq & -\min(2\alpha,\nu\kappa^2)\Vert \widehat{\vu_{\delta_n}}(t,\cdot)\Vert^{2}_{L^2}=-\min(2\alpha,\nu\kappa^2)\Vert \vu_{\delta_n}(t,\cdot)\Vert^{2}_{L^2}.
\end{eqnarray}
Now, we substitute inequality (\ref{ineq_Plancherel}) in  (\ref{ineq_der_time_delta_2}) and we get
$$ \frac{d}{dt} \Vert \vu_{\delta_n}(t,\cdot)\Vert^{2}_{L^2}\leq\frac{\Vert \fe\Vert^{2}_{\dot{H}^{-1}}}{\nu} -\min(2\alpha,\nu\kappa^2)\Vert \vu_{\delta_n}(t,\cdot)\Vert^{2}_{L^2}.$$  We set now $\beta=\min(2\alpha,\nu\kappa^2)>0$ and by the Gr\"onwall inequality we have the control:
$$\Vert \vu_{\delta_n}(t,\cdot)\Vert^{2}_{L^2}\leq e^{-\beta t}\Vert \vu_0 \Vert^{2}_{L^2}+c\frac{\Vert \fe \Vert^{2}_{\dot{H}^{-1}}}{\nu \beta} \left( 1-e^{-\beta t}\right),$$ 
for all time $t\in [0,+\infty[$. Now, we will  recover this control in time  for the limit function $\vu$: we regularize in the time variable the quantity $\Vert \vu_{\delta_n}(t,\cdot)\Vert^{2}_{L^2}$ by  a convolution product with a positive   function $w\in\mathcal{C}^{\infty}_{0}([-\eta,\eta])$ (for $\eta>0$) such that $\ds{\int_{\mathbb{R}}w(t)dt=1}$.  In this way, in the previous inequality we have $$\Vert w\ast \vu_{\delta_n}(t,\cdot)\Vert_{L^2}\leq  w\ast \Vert \vu_{\delta_n}(t,\cdot)\Vert^{2}_{L^2}\leq w\ast \left( e^{-\beta t}\Vert \vu_0 \Vert^{2}_{L^2}+c\frac{\Vert \fe \Vert^{2}_{\dot{H}^{-1}}}{\nu \beta} \left( 1-e^{-\beta t}\right)\right).$$ 
Moreover, since $(\vu_{\delta_n})_{n\in\mathbb{N}}$ converges  weakly$-*$ to $\vu$ in $(L^{\infty}_{t})_{loc}(L^{2}_{x})$ then  $w\ast \vu_{\delta_n}(t,\cdot)$ converges weakly$-*$  to $w\ast\vu(t,\cdot)$ in  $L^2(\Rt)$ and then we can write    
\begin{equation*}
\Vert w\ast \vu(t,\cdot)\Vert^{2}_{L^2} \leq  \liminf_{n\longrightarrow +\infty} \Vert w\ast \vu_{\delta_n}(t,\cdot)\Vert^{2}_{L^2}\leq  w\ast \left( e^{-\beta t}\Vert \vu_0 \Vert^{2}_{L^2}+c\frac{\Vert \fe \Vert^{2}_{\dot{H}^{-1}}}{\nu \beta} \left( 1-e^{-\beta t}\right)\right).
\end{equation*} Now, for $t$ a Lebesgue point of the function $t\mapsto \Vert \vu(t,\cdot)\Vert^{2}_{L^2}$  we have  the control in time  (\ref{control_in_time_Th_2}) and we extend this inequality to all time $t\in [0,+\infty[$ by the weak continuity of the function $t\mapsto \Vert \vu(t,\cdot)\Vert^{2}_{L^2}$. \\

Once we haven proven inequality (\ref{control_in_time_Th_2}), let us deduce from it that $U<+\infty$. Indeed, for a given input scale $\ell_0>0$ we can  write 
$$\ds{\frac{1}{T}\int_{0}^{T} \Vert \vu(t,\cdot)\Vert^{2}_{L^2}\frac{dt}{\ell^{3}_{0}}\leq \frac{1}{T}\int_{0}^{T} e^{-\beta t}\Vert \vu_0 \Vert^{2}_{L^2} \frac{dt}{\ell^{3}_{0}} +\frac{1}{T}\int_{0}^{T}c\frac{\Vert \fe \Vert^{2}_{\dot{H}^{-1}}}{\nu \beta} \left( 1-e^{-\beta t}\right) \frac{dt}{\ell^{3}_{0}}},$$ 
but, since $\fe$ is a time independent function then in the second term in the right-hand  side above we have 
$$ \frac{1}{T}\int_{0}^{T}c\frac{\Vert \fe \Vert^{2}_{\dot{H}^{-1}}}{\nu \beta} \left( 1-e^{-\beta t}\right) \frac{dt}{\ell^{3}_{3}}= c \frac{\Vert \fe \Vert_{\dot{H}^{-1}}}{\nu \beta \ell^{3}_{0}} -  c \frac{1}{T} \int_{0}^{T} \frac{e^{-\beta t}}{\nu \beta} \frac{dt}{\ell^{3}_{0}}.$$ 
Then, taking the limit when $T\longrightarrow +\infty$ we finally obtain 
$$ U^2 =\limsup_{T\longrightarrow +\infty}\frac{1}{T}\int_{0}^{T}\Vert \vu(t,\cdot)\Vert^{2}_{L^2}\frac{dt}{\ell^{3}_{0}}\leq c \frac{\Vert \fe\Vert^{2}_{\dot{H}^{-1}}}{\nu \beta \ell^{3}_{0}}<+\infty,$$ which is estimate  (\ref{estim-U}). \finpv
%%%%%%%%%%%%%%%%%%%%%%%%%%%%%%%%%%
\section{The  Grashof number}\label{Secc_GrashofReynolds} 
We prove here  Theorem \ref{Proposition}. Remark that 
estimates (\ref{estimate_reynolds_grashof}) will be a consequence of the following proposition: 
%%%%%%%%%%%%%%%%%%%%%%%%%%%%%%%%%%%%%%
\begin{Proposition}\label{lemme_estimate_U_F_L} Let  $\vu$ a solution of equation (\ref{damped_alpha_kappa_N-S}) with a damping term $-\alpha P_{\kappa}$ such that $\alpha,\kappa$ satisfy (\ref{alphakappa}). Let $U$ be the characteristic velocity given in (\ref{varepsilon,U}), let $F$ be the averaged force given in (\ref{F}) and let $L$ be the fluid characteristic length given in (\ref{LC}). For a fix number $G_0> 0$, there exists two constants $0<\mathfrak{a}_{1,G_0}\leq \mathfrak{a}_{2,G_0}$, which  depend of the number $G_0$ and $1\leq \theta$, such that we have 
\begin{equation}\label{estim:U-F}
\mathfrak{a}_{1,G_0}\frac{U^2}{L} \leq F \leq \mathfrak{a}_{2,G_0} \frac{U^2}{L}.
\end{equation}
\end{Proposition} 
%%%%%%%%%%%%%%%%%%%%%%%%%%%%%%%%%%%%%%
Indeed, since we have $Re=\frac{UL}{\nu}$ and $Gr=\frac{FL^3}{\nu^{2}}$, we multiply each term of the previous inequality by $\ds{\frac{L^3}{\nu^2}}$ in order to get $\ds{\mathfrak{a}_{1,G_0} \frac{U^2 L^2}{\nu^2} \leq \frac{FL^3}{\nu^2} \leq \mathfrak{a}_{2,G_0} \frac{U^2 L^2}{\nu^2}}$ and thus we obtain the estimate 
$$\ds{\mathfrak{a}_{1,G_0} Re^2 \leq Gr \leq \mathfrak{a}_{2,G_0} Re^2},$$
which is the equivalence announced in Theorem \ref{Proposition}.\\[5mm]
%%%%%%%%%%%%%%%%%%%%%%%%%%%%%%%%%%%%%%
\pv We begin with the inequality $\ds{\mathfrak{a}_{1,G_0} \frac{U^2}{L} \leq F}$ and from the energy inequality (\ref{energy_ineq_alpha_model}) satisfied by the solution $\vu$ (with $\alpha=\frac{\nu}{\ell^{2}_{0}}$ and $\kappa =\frac{1}{20 \theta \ell_{0}}$) we have
\begin{equation}\label{estimate_aux_lemme_U,L,F}
\Vert \vu(t,\cdot)\Vert^{2}_{L^2} \leq \Vert \vu_0 \Vert^{2}_{L^2}+2\int_{0}^{t}\int_{\Rt}\vu(s,x)\cdot\fe(x)dxds +\int_{0}^{t} \left( -2\nu\left\Vert \vu(s,\cdot)\right\Vert^{2}_{\dot{H}^1}  - 2\frac{\nu}{\ell^{2}_{0}} \left\Vert P_{\frac{1}{20 \theta \ell_0}}(\vu)(s,\cdot)\right\Vert^{2}_{L^2} \right)ds.
\end{equation}  
Now,  following the same computations performed in the inequality (\ref{ineq_Plancherel})  for the term inside the last integral above we can write   
$$-2\nu\left\Vert \vu(s,\cdot)\right\Vert^{2}_{\dot{H}^1} -2\frac{\nu}{\ell^{2}_{0}} \left\Vert P(\vu)(s,\cdot)\right\Vert^{2}_{L^2}  \leq  -2\min\left(\frac{\nu}{\ell^{2}_{0}}, \frac{\nu}{20^2 \theta^{2} \ell^{2}_{0}} \right) \Vert \vu(s,\cdot)\Vert^{2}_{L^2} =-\frac{\nu}{200 \,\theta^{2}\, \ell^{2}_{0}} \Vert \vu(s,\cdot)\Vert^{2}_{L^2},$$
and  then, getting back to (\ref{estimate_aux_lemme_U,L,F}) we get
$$ \Vert \vu(t,\cdot)\Vert^{2}_{L^2} \leq \Vert \vu_0 \Vert^{2}_{L^2}+2\int_{0}^{t}\int_{\Rt}\vu(s,x)\cdot\fe(x)dxds -\frac{\nu}{200 \,\theta^{2}\, \ell^{2}_{0}}\int_{0}^{t} \Vert \vu(s,\cdot)\Vert^{2}_{L^2}ds,$$ 
hence we have 
\begin{eqnarray*}
\frac{\nu}{200 \,\theta^{2}\, \ell^{2}_{0}}\int_{0}^{t}\Vert \vu(s,\cdot) \Vert^{2}_{L^2}ds+\Vert \vu(t,\cdot)\Vert^{2}_{L^2}&\leq & \Vert \vu_0 \Vert^{2}_{L^2}+2\int_{0}^{t}\int_{\Rt}\vu(s,x)\cdot\fe(x)dxds\\ 
\frac{\nu}{200 \,\theta^{2}\, \ell^{2}_{0}}\int_{0}^{t}\Vert \vu(s,\cdot) \Vert^{2}_{L^2}ds&\leq&  \Vert \vu_0 \Vert^{2}_{L^2}+2\int_{0}^{t}\int_{\Rt}\vu(s,x)\cdot\fe(x)dxds.
\end{eqnarray*}
%Revisar con cuidado
%%%%%%%%%%%%%%%%%%%%%%%%%%%%%%%%%%%%%%%%%%%%%
Now, in the second term in the right side, we apply the Cauchy-Schwarz inequality and the Young inequalities to obtain
\begin{eqnarray*}
\frac{\nu}{200 \,\theta^{2}\, \ell^{2}_{0}}\int_{0}^{t}\Vert \vu(s,\cdot) \Vert^{2}_{L^2}ds&\leq  & \Vert \vu_0 \Vert^{2}_{L^2}+2\int_{0}^{t}\Vert \fe \Vert_{L^2}\Vert \vu(s,\cdot)\Vert_{L^2} ds\\
& \leq & \Vert \vu_0 \Vert^{2}_{L^2}+\frac{400 \, \theta^2  \ell^{2}_{0}}{\nu }t\Vert \fe \Vert^{2}_{L^2}+\frac{\nu}{400 \, \theta^2 \ell^{2}_{0}}\int_{0}^{t}\Vert \vu(s,\cdot)\Vert^{2}_{L^2}ds,
\end{eqnarray*} 
from which we deduce the estimate	
$$ \frac{\nu}{400 \, \theta^2 \ell^{2}_{0}} \int_{0}^{t} \Vert \vu(s,\cdot)\Vert^{2}_{L^2}ds \leq \Vert \vu_0 \Vert^{2}_{L^2} +\frac{400\, \theta^{2} \ell^{2}_{0}}{\nu} t \Vert \fe \Vert^{2}_{L^2}. $$
Now, dividing by $t$ and $\ell^{3}_{0}$ each term above and taking the limit $\underset{t\longrightarrow+\infty}{\limsup}$ we get 
$$ \frac{\nu}{400 \, \theta^{2} \ell^{2}_{0}} \limsup_{t\longrightarrow +\infty}\frac{1}{t}\int_{0}^{t}\Vert \vu(s,\cdot) \Vert^{2}_{L^2}\frac{ds}{\ell^{3}_{0}} \leq  \frac{400 \, \theta^{2}\ell^{2}_{0}}{\nu}\frac{\Vert \fe \Vert^{2}_{L^2}}{\ell^{3}_{0}}.$$ 
Finally, using the definitions of $U$ and $F$ given in formulas (\ref{varepsilon,U}) and (\ref{F}) respectively one obtains, after a division by $L$
$$ \frac{1}{400^2 \theta^4} \frac{U^2}{L}\leq \frac{\,\ell^{4}_{0}}{\nu^2}\frac{F^2}{L}.$$
Observe now that the right-hand side above can be rewritten in the form
$\frac{\ell^{4}_{0}}{\nu^2}\frac{F^2}{L}=\left( \left(\frac{\ell^{4}_{0}}{L^4}\right)\left(\frac{F L^3}{\nu^2}\right)\right)F$, but since $\frac{\ell^{4}_{0}}{L^4}= \gamma^4$ by (\ref{LC}) and $Gr=\frac{FL^3}{\nu^2}$ by (\ref{Grashof}) and $\gamma^4 Gr=\frac{G_0}{c_0}$ by (\ref{Grashof_zeo}), we actually have $\frac{\ell^{4}_{0}}{\nu^2}\frac{F^2}{L}=\frac{G_0}{c_0} F$, and thus we obtain from the previous inequality:
$$ \frac{c_0}{400^2 \, \theta^4 G_0} \frac{U^2}{L}\leq F.$$
It remains to set the constant  
\begin{equation}\label{C_1}
0<\mathfrak{a}_{1,G_0}\leq \frac{c_0}{400^2 \, \theta^4 G_0},
\end{equation} 
to obtain the desired estimate $\mathfrak{a}_{1,G_0} \frac{U^2}{L} \leq F$. \\

We continue now with the upper estimate $\ds{F \leq \mathfrak{a}_{2,G_0} \frac{U^2}{L}}$: observe that since $\vu\in L^{\infty}_{t}(L^{2}_{x})\cap (L^{2}_{t})_{loc}(\dot{H}^{1}_{x})$, we have $\partial_t \vu \in (L^{2}_{t})_{loc}(H^{-\frac{3}{2}}_{x})$,  $\P((\vu\cdot \vec{\nabla}) \vu)\in (L^{2}_{t})_{loc}(H^{-\frac{3}{2}}_{x})$,  $\Delta \vu \in (L^{2}_{t})_{loc}(H^{-1}_{x})$ and  $P_{\frac{1}{20\, \theta \ell_0}}(\vu)\in L^{\infty}_{t}(L^{2}_{x})$. On the other hand,  since $\widehat{\fe}$ is localized at the frequencies $\frac{1}{10\theta\ell_0}\leq \vert \xi \vert \leq \frac{1}{\theta\ell_0}$, then $\fe$ belongs to all  Sobolev spaces $H^{s}(\Rt)$ with $s\in \mathbb{R}$, thus integrating in the space variable we have the identity
\begin{eqnarray*}
\int_{\Rt} \partial_t\vu(t,x)\cdot \fe(x)dx&= &\int_{\Rt}\nu\Delta \vu(t,x)\cdot\fe(x)dx-\int_{\Rt}\P((\vu\cdot\vec{\nabla})\vu)(t,x)\cdot \fe(x)dx+\Vert \fe \Vert^{2}_{L^2}\\
& & -\frac{\nu}{\ell^{2}_{0}} \int_{\Rt} P_{\frac{1}{20\, \theta \ell_0}}(\vu)(t,x)\cdot \fe(x)dx,
\end{eqnarray*} 
which can be rewritten in the following form
\begin{eqnarray}\label{eq_aux_main_theorem} \nonumber 
\Vert \fe \Vert^{2}_{L^2} &=& \int_{\Rt} \partial_t\vu(t,x)\cdot \fe(x)dx -\int_{\Rt}\nu\Delta \vu(t,x)\cdot\fe(x)dx+\int_{\Rt}\P((\vu\cdot\vec{\nabla})\vu)(t,x)\cdot \fe(x)dx\\
& & +\frac{\nu}{\ell^{2}_{0}} \int_{\Rt} P_{\frac{1}{20\, \theta \ell_0}}(\vu)(t,x)\cdot \fe(x)dx.
\end{eqnarray}
Note that now, since $\fe$ is stationary, for the first term in the right-hand side we have 
$$\ds{\int_{\Rt} \partial_t\vu(t,x)\cdot \fe(x)dx= \partial_t\int_{\Rt}\vu(t,x)\cdot \fe(x)dx}.$$ 
For the second term in the right-hand side of (\ref{eq_aux_main_theorem}), integrating by parts and applying the Cauchy-Schwarz inequality we get 
$$-\int_{\Rt}\nu\Delta \vu(t,x)\cdot\fe(x)dx =-\nu \int_{\Rt}\vu(t,x)\cdot\Delta\fe(x)dx\leq \nu  \Vert\vu(t,\cdot)\Vert_{L^2}\Vert \Delta \fe\Vert_{L^2}.$$ 
For the third term in the right-hand side of (\ref{eq_aux_main_theorem}), since $\fe$ is a divergence-free function, by the properties of the Leray projector, with an integration by parts and by the H\"older inequalities we obtain
\begin{eqnarray*}
\int_{\Rt}\P((\vu\cdot\vec{\nabla})\vu)(t,x)\cdot \fe(x)dx&=& \int_{\Rt}(\vu\cdot\vec{\nabla})\vu(t,x)\cdot \fe(x)dx=-\sum_{i,j=1}^{3}\int_{\Rt}u_{i}(t,x)u_{j}(t,x)\partial_{j}f_i(x)dx\\
& \leq & \Vert \vu(t,\cdot)\Vert^{2}_{L^2}\Vert \vec{\nabla} \otimes\fe\Vert_{L^{\infty}}.
\end{eqnarray*}
Finally, for the fourth term in the right-hand side of (\ref{eq_aux_main_theorem}), recall that $\kappa=\frac{1}{20\, \theta \ell_0}$ and if $1\leq \theta$ then we have 
\begin{eqnarray*}
\int_{\Rt}P_{\frac{1}{20 \, \theta\,  \ell_0}}(\vu)(t,x)\cdot \fe(x)\, dx &=&\int_{\Rt}\widehat{P_{\frac{1}{20 \,\theta   \ell_0}}(\vu)}(t,\xi)\cdot \widehat{\fe}(\xi)\,d\xi\\
&=&\int_{\Rt}\left(\mathds{1}_{\vert \xi \vert<\frac{1}{20\, \theta  \ell_0}}(\xi)\widehat{\vu}(t,\xi)\right)\cdot\left( \mathds{1}_{\frac{1}{10\theta\ell_0}\leq\vert \xi \vert\leq\frac{1}{\theta\ell_0}}(\xi)\widehat{\fe}(\xi)\right)d\xi=0.
\end{eqnarray*} 
Then, gathering all these previous remarks, the expression (\ref{eq_aux_main_theorem}) becomes
$$ \Vert \fe \Vert^{2}_{L^2}\leq  \partial_t \int_{\Rt} \vu(t,x)\cdot \fe(x)dx+ \Vert \vu(t,\cdot)\Vert^{2}_{L^2}\Vert \vec{\nabla} \otimes\fe\Vert_{L^{\infty}}+\nu\Vert\vu(t,\cdot)\Vert_{L^2}\Vert \Delta \fe\Vert_{L^2}.$$
Now, for $T>0$ we take the average $\ds{\frac{1}{T}\int_{0}^{T}(\cdot)dt}$ and applying the Cauchy-Schwarz inequality (in time variable) in the last term above and using the fact that $\fe$ is stationary we have 
\begin{eqnarray*}
\Vert \fe \Vert_{L^2}^2 &\leq &\frac{1}{T} \left( \int_{\Rt} \vu(T,x)\cdot \fe(x)dx- \vu_{0}(x)\cdot \fe(x)dx \right)+\left( \frac{1}{T} \int_{0}^{T} \Vert \vu(t,\cdot)\Vert^{2}_{L^2} dt \right) \Vert \vec{\nabla} \otimes \fe \Vert_{L^{\infty}}\\
& & +\nu \left( \frac{1}{T} \int_{0}^{T} \Vert \vu(t,\cdot)\Vert_{L^2}^{2}dt \right)^{\frac{1}{2}}\Vert \Delta \fe \Vert_{L^2},
\end{eqnarray*} 
and thus, taking the limit $\ds{\limsup_{T\longrightarrow+\infty}}$ and by the definition of the characteristic velocity $U$ given in (\ref{varepsilon,U}), one has
\begin{equation}\label{eq_3_aux_main_theorem}
\Vert \fe \Vert^{2}_{L^2} \leq \limsup_{T\longrightarrow +\infty}\frac{1}{T}\left( \int_{\Rt}\vu(T,x)\cdot\fe(x)-\vu_0(x)\cdot\fe(x)dx\right)+\ell^{3}_{0}U^{2}\Vert \vec{\nabla} \otimes \fe \Vert_{L^{\infty}} +\nu \ell^{3/2}_{0}U \Vert \Delta \fe \Vert_{L^2}.
\end{equation}
We will prove now that the first term in the right-hand side above is actually null. Indeed, applying the Cauchy-Schwarz inequality we get 
$$\frac{1}{T}\int_{\Rt}\vu(T,x)\cdot\fe(x)-\vu_0(x)\cdot\fe(x) dx\leq \frac{1}{T}\big( \Vert \vu(T,\cdot)\Vert_{L^2} +\Vert \vu_0\Vert_{L^2}\big)\Vert \fe \Vert_{L^2},$$ 
and since the velocity $\vu$ satisfies the estimate (\ref{control_in_time_Th_2}), we have the control
$$ \frac{1}{T}\big( \Vert \vu(T,\cdot)\Vert_{L^2} +\Vert \vu_0\Vert_{L^2}\big)\Vert \fe \Vert_{L^2}\leq  \frac{1}{T}\left(e^{-\beta T}\Vert \vu_0 \Vert^{2}_{L^2}+c\frac{\Vert \fe \Vert^{2}_{\dot{H}^{-1}}}{\nu \beta} \left( 1-e^{-\beta T}\right)+\Vert \vu_0\Vert_{L^2}\right)\Vert \fe \Vert_{L^2},$$
from which we easily deduce the that the first term of (\ref{eq_3_aux_main_theorem}) is null. Thus, dividing (\ref{eq_3_aux_main_theorem})  by $\ell^{3}_{0}$ and by the definition of $F$ given in (\ref{F}) we obtain the inequality $\ds{F^2\leq U^2\Vert \vec{\nabla} \otimes \fe \Vert_{L^{\infty}}+\nu U \frac{\Vert \Delta \fe \Vert_{L^2}}{\ell^{\frac{3}{2}}_{0}}},$ which is equivalent to 
\begin{equation}\label{averaged_ineq_0}
F\leq U^2 \frac{\Vert \vec{\nabla} \otimes \fe \Vert_{L^{\infty}}}{F}+\nu U \frac{\Vert \Delta \fe \Vert_{L^2}}{\ell^{\frac{3}{2}}_{0}F}.
\end{equation}
Observe now that since  $\widehat{\fe}$ is localized at the frequencies $\frac{1}{10\theta\ell_0}\leq \vert \xi \vert \leq\frac{1}{\theta\ell_0}$ then by the Bernstein inequalities one has the equivalence $\frac{\Vert \vec{\nabla} \otimes \fe \Vert_{L^{\infty}}}{F}\simeq \frac{1}{\theta\ell_0}\frac{\Vert \fe \Vert_{L^{\infty}}}{F}\lesssim  \frac{1}{\ell_0}\frac{\Vert \fe \Vert_{L^{\infty}}}{F}$. Moreover using the formulas (\ref{LC}), (\ref{gamma}) and (\ref{F}) we obtain $\frac{\Vert \vec{\nabla} \otimes \fe \Vert_{L^{\infty}}}{F} \leq \frac{c_{1}}{L}$. 
On the other hand, always by the Bernstein inequalities, we have  $\Vert \Delta \fe \Vert_{L^2}\simeq \left( \frac{1}{\theta\ell_0}\right)^2\Vert \fe \Vert_{L^2}$ and since $\ell^{\frac{3}{2}}_{0}F=\Vert \fe \Vert_{L^2}$ we can write $\frac{\Vert \Delta \fe \Vert_{L^2}}{\ell^{\frac{3}{2}}_{0}F} \simeq \left(\frac{1}{\theta\ell_0}\right)^2=\frac{1}{\theta^{2}\gamma^2L^2}$. With this two facts in mind, we can rewrite inequality (\ref{averaged_ineq_0}) in the following manner:
$$F\leq  c_{1}\frac{U^2}{L}+\frac{c_{2}}{\theta^{2}}\frac{\nu U}{\gamma^2L^2},$$
where $c_{1}$ and $c_{2}$ are the constants that came from the Bernstein inequalities. Thus, if we denote $c_{3}=\frac{c_{2}}{\theta^{2}}$ and since $Gr=\frac{FL^3}{\nu^2}$ and $Gr\gamma^4= \frac{G_0}{c_0} $, we have
\begin{equation}\label{averaged_ineq_01}
0\leq c_{1}\frac{U^2}{FL} +c_{3}\frac{U}{(FL)^{\frac{1}{2}}}\frac{(c_0)^{\frac{1}{2}}}{(G_0)^{\frac{1}{2}}}-1.
\end{equation}
Now we fix the technical parameter $1\leq \theta$ large enough such that 
\begin{equation}\label{Technical}
c_{3}=\frac{c_{2}}{\theta^{2}}<\frac{1}{2},
\end{equation}
and this upper bound will be useful in the sequel. We continue the study of estimate (\ref{averaged_ineq_01}) and we define the variable  $x=\frac{U}{(FL)^{\frac{1}{2}}}$, thus solving the equation $\ds{ 0\leq c_{1}x^2 +\frac{c_3 (c_0)^{\frac{1}{2}}}{G^{\frac{1}{2}}_{0}}x-1}$ we obtain the constraint $\ds{ \frac{1}{2 c_{1}}\left(-\frac{c_3 \sqrt{c_0}}{\sqrt{G_0}}+\sqrt{\frac{{c_3}^2 c_0}{G_0}+ 4c_{1}}\right) \leq x}$, and getting back to the initial variables we have
$$ \frac{1}{2 c_{1}}\left(-\frac{c_3 \sqrt{c_0}}{\sqrt{G_0}}+\sqrt{\frac{{c_3}^2 c_0}{G_0}+ 4c_{1}}\right)   \sqrt{F}\leq \frac{U}{\sqrt{L}}.$$
It only remains to set the constant 
\begin{equation}\label{C_2}
\ds{\mathfrak{a}_{2,G_0}= \left(\frac{1}{2 c_{1}}\left(-\frac{c_3 \sqrt{c_0}}{\sqrt{G_0}}+\sqrt{\frac{{c_3}^2 c_0}{G_0}+ 4c_{1}}\right) \right)^{-2}},
\end{equation}
and then we obtain the upper estimate $F\leq \mathfrak{a}_{2,G_0} \frac{U^2}{L}$.
To conclude, we need to verify the compatibility of the constants $\mathfrak{a}_{1,G_0}$ and $\mathfrak{a}_{2,G_0}$ given in (\ref{C_1}) and (\ref{C_2}), \emph{i.e.} we must check that $\mathfrak{a}_{1,G_0}\leq \mathfrak{a}_{2,G_0}$ and this condition is equivalent to 
$$\mathfrak{a}_{1,G_0}\leq \frac{c_0}{400^2 \, \theta^4 G_0} \leq  \left(\frac{1}{2 c_{1}}\left(-\frac{c_3 \sqrt{c_0}}{\sqrt{G_0}}+\sqrt{\frac{{c_3}^2 c_0}{G_0}+ 4c_{1}}\right) \right)^{-2}=\mathfrak{a}_{2,G_0},$$
which is satisfied as long as $\frac{c_{0}(1+\sqrt{c_3})^{2}}{c_{1}}<G_{0}$. 
Proposition \ref{lemme_estimate_U_F_L} is now proven.   \finpv	
%%%%%%%%%%%%%%%%%%%%%%%%%%%%%%%%%%%%%%
\section{Proof of Theorem \ref{Main_Result}}\label{Sec_Kolmogorov} 
We decompose the proof of the Kolmogorov dissipation law in two steps and we start with the lower estimate 
$$\ds{\mathfrak{b}_{1,G_0}\frac{U^3}{L}\leq \varepsilon}.$$
Recall that in the previous Proposition \ref{lemme_estimate_U_F_L} we proved the inequality $\mathfrak{a}_{1,G_0}\frac{U^2}{L}\leq  F$, or equivalently 
\begin{equation}\label{EstimationBase1}
\ds{\mathfrak{a}_{1,G_0}\frac{U^3}{L}\leq FU},
\end{equation}
and now we will prove the following estimate \begin{equation}\label{ineq_FU_varep}
\ds{FU \leq  (20\, \theta)^{2} \varepsilon}. 
\end{equation}
This estimate will be attained as long as the  Grashof number $Gr$ is large enough and satisfies the condition $\ds{\frac{4  \mathfrak{a}_{2,G_0} }{c^{2}_{0}} \frac{G^{2}_{0}}{\gamma^4} \leq Gr}$,  but recall that  we have the identity $\ds{Gr=\frac{G_0}{c_0 \gamma^4}}$ and thus, this condition is verified when the number $G_0$ satisfies:
\begin{equation}\label{Eq_Condition_Grashof}
\frac{4 \mathfrak{a}_{2,G_0}}{c_0} G_0 \leq 1.
\end{equation}  Recall also that the constant  $\mathfrak{a}_{2,G_0}$ depends on the number $G_0$ and is given in (\ref{C_2}) and thus, (\ref{Eq_Condition_Grashof})  is equivalent to the constraint $4\frac{c_{1}G_{0}}{c_{0}c_{3}}\leq \left(\sqrt{1+\frac{4}{c_{3}}\frac{c_{1}G_{0}}{c_{0}c_{3}}}-1 \right)$ which is valid as long as $G_{0}$ is  small and $c_{3}$ satisfies (\ref{Technical}). \\

We assume from now on the control (\ref{Eq_Condition_Grashof}) and therefore  we have $\ds{\frac{4  \mathfrak{a}_{2,G_0} }{c^{2}_{0}} \frac{G^{2}_{0}}{\gamma^4} \leq Gr}$ from which we can write  we can write $\frac{4}{c^{2}_{0}} \frac{G^{2}_{0}}{\gamma^4}  \leq \frac{1}{\mathfrak{a}_{2,G_0}} Gr$ and since by Theorem \ref{Proposition} we have the inequality $\ds{ Gr \leq  \mathfrak{a}_{2,G_0} \,Re^2}$, we can deduce from these two facts the control $\ds{\frac{2}{c_0} \frac{G_0}{\gamma^2}  \leq Re=\frac{UL}{\nu}}$. Moreover by (\ref{gamma}) and (\ref{F}) and (\ref{Grashof_zeo}) we can write $\ds{\frac{2}{\gamma}\frac{F \ell^{3}_{0}}{\nu^2}\leq \frac{UL}{\nu}}$ and using the definition of the characteristic length $L=\frac{\ell_0}{\gamma}$ we actually obtain  $\ds{\frac{2F \ell^{2}_{0}}{\nu}\leq U}$ which is equivalent to 
$$2FU\leq \frac{\nu}{\ell^{2}_{0}}U^2.$$
Thus, in order to prove (\ref{ineq_FU_varep}) we will show that we have 
\begin{equation}\label{averaged_ineq_1}
\frac{\nu}{\ell^{2}_{0}}U^2 \leq FU+ (20\, \theta)^{2} \varepsilon.   
\end{equation} 
To this end, we use the energy inequality (\ref{energy_ineq_alpha_model}) and the definition of the operator $-\alpha P_{\kappa}$ where by hypothesis we have $\alpha= \frac{\nu}{\ell^{2}_{0}}$ and $\kappa= \frac{1}{20\, \theta\ell_0}$ (see (\ref{alphakappa})):
\begin{eqnarray*}
\Vert \vu(t,\cdot)\Vert^{2}_{L^2}&\leq &-2\nu\int_{0}^{t} \Vert \vec{\nabla} \otimes \vu(s,\cdot)\Vert^{2}_{L^2}ds + \Vert \vu_0 \Vert^{2}_{L^2}+2\int_{0}^{t}\int_{\Rt}\vu(s,x)\cdot\fe(x)dxds\\
& &  -2\frac{\nu}{\ell^{2}_{0}} \int_{0}^{t}\Vert P_{\frac{1}{20\, \theta \ell_0}}(\vu)(s,\cdot) \Vert^{2}_{L^2}ds \\ 
&\leq & \Vert \vu_0 \Vert^{2}_{L^2}+2\int_{0}^{t}\int_{\Rt}\vu(s,x)\cdot\fe(x)dxds  -2\frac{\nu}{\ell^{2}_{0}} \int_{0}^{t}\Vert \mathds{1}_{\vert \xi\vert<\frac{1}{20\, \theta \ell_0}}\widehat{\vu}(s,\xi)\Vert^{2}_{L^2}ds.
\end{eqnarray*}
Introducing useful information we have
\begin{eqnarray*}
&\leq & \Vert \vu_0 \Vert^{2}_{L^2}+2\int_{0}^{t}\int_{\Rt}\vu(s,x)\cdot\fe(x)dxds  -2\frac{\nu}{\ell^{2}_{0}} \int_{0}^{t}\Vert \mathds{1}_{\vert \xi\vert<\frac{1}{20\, \theta \ell_0}}\widehat{\vu}(s,\xi)\Vert^{2}_{L^2}ds\nonumber\\
& &-2\frac{\nu}{\ell^{2}_{0}} \int_{0}^{t}\Vert \mathds{1}_{\vert \xi\vert\geq \frac{1}{20\, \theta \ell_0}}\widehat{\vu}(s,\xi)\Vert^{2}_{L^2}ds+2\frac{\nu}{\ell^{2}_{0}} \int_{0}^{t}\Vert \mathds{1}_{\vert \xi\vert\geq \frac{1}{200\ell_0}}\widehat{\vu}(s,\xi)\Vert^{2}_{L^2}ds \nonumber\\ 
&\leq & \Vert \vu_0 \Vert^{2}_{L^2}+2\int_{0}^{t}\int_{\Rt}\vu(s,x)\cdot\fe(x)dxds  -2\frac{\nu}{\ell^{2}_{0}}\int_{0}^{t} \Vert \vu(s,\cdot)\Vert^{2}_{L^2}ds+2\frac{\nu}{\ell^{2}_{0}} \int_{0}^{t}\Vert \mathds{1}_{\vert \xi\vert\geq \frac{1}{20\, \theta \ell_0}}\widehat{\vu}(s,\xi)\Vert^{2}_{L^2}ds ,
\end{eqnarray*}
and we obtain
\begin{eqnarray*}
2\frac{\nu}{\ell^{2}_{0}} \int_{0}^{t}\Vert \vu(s,\cdot) \Vert^{2}_{L^2}ds +\Vert \vu(t,\cdot)\Vert^{2}_{L^2} &\leq & \Vert \vu_0 \Vert^{2}_{L^2}+2\int_{0}^{t}\int_{\Rt}\vu(s,x)\cdot\fe(x)dxds\nonumber\\
& &+2\frac{\nu}{\ell^{2}_{0}} \int_{0}^{t}\left\Vert \mathds{1}_{\vert \xi \vert\geq \frac{1}{20\, \theta \ell_0}}\widehat{\vu}(s,\cdot)\right\Vert^{2}_{L^2}ds
\end{eqnarray*}
\begin{eqnarray}
2\frac{\nu}{\ell^{2}_{0}} \int_{0}^{t}\Vert \vu(s,\cdot) \Vert^{2}_{L^2}ds \leq  \Vert \vu_0 \Vert^{2}_{L^2}+2\int_{0}^{t}\int_{\Rt}\vu(s,x)\cdot\fe(x)dxds
+2\frac{\nu}{\ell^{2}_{0}} \int_{0}^{t}\left\Vert \mathds{1}_{\vert \xi \vert\geq \frac{1}{200\ell_0}}\widehat{\vu}(s,\cdot)\right\Vert^{2}_{L^2}ds.\label{ineq_tech_Th_auxiliar}
\end{eqnarray}
We study now the last term above and we write
\begin{eqnarray*}
2\frac{\nu}{\ell^{2}_{0}} \int_{0}^{t}\left\Vert \mathds{1}_{\vert \xi \vert\geq \frac{1}{20\, \theta \ell_0}}\widehat{\vu}(s,\cdot)\right\Vert^{2}_{L^2}ds&=& 2\times(20\, \theta)^{2}\nu \int_{0}^{t}\left\Vert  \left( \frac{1}{20\, \theta \ell_{0}}\right) \mathds{1}_{\vert \xi \vert\geq \frac{1}{20\, \theta \ell_0}}\widehat{\vu}(s,\cdot)\right\Vert^{2}_{L^2}ds\\
& \leq & 2\times(20\, \theta)^{2}\nu \int_{0}^{t}\left\Vert \vert \xi \vert \widehat{\vu}(s,\cdot)\right\Vert^{2}_{L^2} ds=2\times(20\, \theta)^{2} \nu \int_{0}^{t} \Vert \vu(s,\cdot)\Vert^{2}_{\dot{H}^1}ds,
\end{eqnarray*} 
and thus, coming back to (\ref{ineq_tech_Th_auxiliar}) we have
\begin{eqnarray*}
2\frac{\nu}{\ell^{2}_{0}}\int_{0}^{t}\Vert \vu(s,\cdot) \Vert^{2}_{L^2}ds& \leq & \Vert \vu_0 \Vert^{2}_{L^2}+2\int_{0}^{t}\int_{\Rt}\vu(s,x)\cdot \fe(x)  dxds +2\times (20\, \theta)^{2}\nu\int_{0}^{t} \Vert\vu (s,\cdot) \Vert^{2}_{\dot{H}^1}ds\\
&\leq &  \Vert \vu_0 \Vert^{2}_{L^2} +2\int_{0}^{t}\Vert \fe \Vert_{L^2} \Vert \vu(s,\cdot)\Vert_{L^2}ds +2\times(20\, \theta)^{2}\nu\int_{0}^{t} \Vert \vu (s,\cdot) \Vert^{2}_{\dot{H}^1}ds\\
&\leq & \Vert \vu_0 \Vert^{2}_{L^2} +2\Vert \fe \Vert_{L^2} \left( \int_{0}^{t}\Vert \vu(s,\cdot)\Vert^{2}_{L^2}ds\right)^{\frac{1}{2}}t^{\frac{1}{2}}+ 2\times (20\, \theta)^{2}\nu\int_{0}^{t} \Vert \vu (s,\cdot) \Vert^{2}_{\dot{H}^1}ds,
\end{eqnarray*} 
dividing the previous expression by $\ell^{3}_{0}$ and $t$ and taking the limit $\underset{t\longrightarrow +\infty}{\limsup}$, with the definition of $\varepsilon$, $U$ and $F$ given in (\ref{varepsilon,U}) and (\ref{F}) respectively, we have  
$$ \frac{\nu}{\ell^{2}_{0}}U^{2} \leq FU+(20\, \theta)^{2} \varepsilon,$$
which is exactly the estimate given in (\ref{averaged_ineq_1}): it only remains to set $\mathfrak{b}_{1,G_0}=\frac{\mathfrak{a}_{1,G_0}}{(20\, \theta)^{2}}$.\\

We study now the upper estimate $\ds{\varepsilon \leq \mathfrak{b}_{2,G_0} \frac{U^3}{L}}$. Note that by Proposition \ref{lemme_estimate_U_F_L} we have $\ds{F \leq \mathfrak{a}_{2,G_0} \frac{U^2}{L}}$ which is equivalent to $FU \leq \mathfrak{a}_{2,G_0}\frac{U^3}{L}$ and then it only remains to prove  $\ds{\varepsilon \leq FU}$ and to set $\mathfrak{b}_{2,G_0}=\mathfrak{a}_{2,G_0}$. Using the energy inequality  (\ref{energy_ineq_alpha_model}) we can write: 
\begin{eqnarray*}
\Vert \vu(t,\cdot)\Vert^{2}_{L^2}+2\nu\int_{0}^{t} \Vert \vu(s,\cdot)\Vert^{2}_{\dot{H}^1}ds&\leq &\Vert \vu_0 \Vert^{2}_{L^2}+2\int_{0}^{t} \langle \fe, \vu(s,\cdot)\rangle_{\dot{H}^{-1}\times \dot{H}^1} ds -2\alpha \int_{0}^{t}\Vert P_{\frac{1}{20\, \theta \ell_0}}(\vu)(s,\cdot) \Vert^{2}_{L^2}ds\\
2\nu\int_{0}^{t} \Vert \vu(s,\cdot)\Vert^{2}_{\dot{H}^1}ds&\leq&\Vert \vu_0 \Vert^{2}_{L^2}+ 2\int_{0}^{t}\int_{\Rt}\fe(x)\cdot\vu(s,x)dxds,
\end{eqnarray*}
applying the Cauchy-Schwarz inequality in the last term above (first in spatial variable and then in time variable) we obtain
$$2\nu\int_{0}^{t} \Vert \vu(s,\cdot)\Vert^{2}_{\dot{H}^1}ds \leq  \Vert \vu_0 \Vert^{2}_{L^2} +2t^{\frac{1}{2}}\left( \int_{0}^{t}\Vert \vu(s,\cdot)\Vert^{2}_{L^2}ds\right)^{\frac{1}{2}}\Vert \fe \Vert_{L^2},$$ 
and dividing by  $\ell^{3}_{0}$ and $T$, and taking the limit when $t \longrightarrow +\infty$, we get by the definition of $\varepsilon$, $U$ and $F$ (see  (\ref{varepsilon,U}) and (\ref{F})) the whished inequality $\ds{\varepsilon \leq FU}$.  Theorem \ref{Main_Result} is now proven.  \finpv

%%%%%%%%%%%%%%%%%%%%%%%%%%%%%%%%%%%%%%%%%%%%%%%%%%%%%%%%%%%%%%%%%%%%%%%%%%%
\section{A non-turbulent model}\label{Sec_non_turbulent}
\subsection{Proof of Theorem \ref{Th-non-turb}} 
Recall that by Theorem \ref{Main_Result} we have the inequalities $\ds{\mathfrak{b}_{1,G_{0}} \frac{U^3}{L} \leq \varepsilon \leq  \mathfrak{b}_{2,G_{0}} \frac{U^3}{L}}$ and using the definition of the Taylor scale given in (\ref{Taylor_scale}) we have 
$$\ds{\left( \frac{1}{\mathfrak{b}_{2,G_{0}}} \frac{\nu }{U L} \right)^{\frac{1}{2}} L \leq \ell_T \leq \left( \frac{1}{\mathfrak{b}_{1,G_{0}}} \frac{\nu }{U L} \right)^{\frac{1}{2}} L}.$$
Now, since $\ds{Re=\frac{UL}{\nu}}$ and since by Theorem \ref{Proposition} we have $\ds{\mathfrak{a}_{1,G_{0}}^{\frac{1}{4}} Re^{\frac{1}{2}} \leq G^{\frac{1}{4}}_{r} \leq \mathfrak{a}_{2,G_{0}}^{\frac{1}{4}} Re^{\frac{1}{2}}}$ we can write
\begin{equation}\label{estim-Taylor2}
\left( \frac{\mathfrak{a}_{1,G_{0}}^{\frac{1}{4}}}{\mathfrak{b}_{2,G_{0}}^{\frac{1}{2}}} \right)   \frac{L}{G^{\frac{1}{4}}_{r}} \leq \ell_T  \leq \left( \frac{\mathfrak{a}_{2,G_{0}}^{\frac{1}{4}}}{\mathfrak{b}_{1,G_{0}}^{\frac{1}{2}}} \right)   \frac{L}{G^{\frac{1}{4}}_{r}},
\end{equation} 
but, recall that we have $L=\frac{\ell_0}{\gamma}$ (see (\ref{LC})) and $\ds{Gr=\frac{G_0}{c_0 \gamma^4}}$ (see (\ref{estimate_Grashof_G0})), we can write  $\ds{\frac{L}{G^{\frac{1}{4}}_{r}}= \ell_0 \left( \frac{c^{\frac{1}{4}}_{0}}{G^{\frac{1}{4}}_{0}} \right)}$, and thus coming back to the expression (\ref{estim-Taylor2}) we finally have 
$$\ds{ \left(  \frac{\mathfrak{a}_{1,G_{0}}^{\frac{1}{4}}}{\mathfrak{b}_{2,G_{0}}^{\frac{1}{2}}} \frac{c^{\frac{1}{4}}_{0}}{G^{\frac{1}{4}}_{0}} \right) \ell_0 \leq \ell_T \leq \left( \frac{\mathfrak{b}_{2,G_{0}}^{\frac{1}{4}}}{\mathfrak{b}_{1,G_{0}}^{\frac{1}{2}}} \frac{c^{\frac{1}{4}}_{0}}{G^{\frac{1}{4}}_{0}}\right) \ell_0},$$ 
and Theorem  \ref{Th-non-turb} is proven. \finpv 

%%%%%%%%%%%%%%%%%%%%%%%%%%%%%%%%%%%%%%
\subsection{Example with a particular external force}\label{Sec:exemple_force}
Remark that the number $G_0$ is defined by means of the force $\fe$ in expression (\ref{Grashof_zeo}) and the fact that we  fixed this number implies a control on the  amplitude of the  force. Indeed, by expression (\ref{Grashof_zeo}) we write $\Vert \fe \Vert_{L^{\infty}}= \frac{G_0 \nu^2}{\ell^{3}_{0}}$ and since $G_0$ is  fixed then we have the equivalence  $\Vert \fe \Vert_{L^{\infty}}\simeq \frac{\nu^2}{\ell^{3}_{0}}$. Moreover, remark also  that we have $Gr\gg1$ if $\gamma\ll1$ which, by (\ref{gamma}) and the fact that the quantity $\Vert \fe \Vert_{L^{\infty}}$ is now fixed as above,  is equivalent to the condition $  \Vert \fe \Vert_{L^{2}} \gg 1$.  \\

In this section  we construct an example   of a  force $\fe$ with the property $\Vert \fe \Vert_{L^2}\gg 1$ and which its amplitude verifies the equivalence above.  We will use for this a wavelet. 
\begin{Definition}\label{Ondelette_base} 
Let  $\vec{\phi}=(\phi_1, \phi_2, \phi_3)$ be a vector field  in the Schwartz class such that:
\begin{enumerate}
\item[1)] $\vec{\phi}$ is a  divergence-free vector,
\item[2)] for $1\leq \theta$, we have   $supp\,(\widehat{\phi_i})\subset \lbrace\xi\in \Rt: \frac{1}{10\theta} \leq \vert \xi \vert\leq \frac{1}{\theta} \rbrace,$ for all $i=1,2,3$,
\item[3)] for all $\vec{\psi}$ vector field in the Schwartz class which verifies the property $2)$ above we have 
$$\displaystyle{\int_{\Rt}\vec{\phi}(x-k)\cdot \vec{\psi}(x-m)dx}=\delta_{k,m},$$ for all $k,m\in \Zt,$ and where $\delta_{k,m}$ is the Kronecker's delta function. 
\end{enumerate}
\end{Definition} 
We state now a well-known property of wavelets. 
\begin{Lemme}\label{prop_ondelette_base} Let $\vec{\phi}=(\phi_1,\phi_2,\phi_2)$ be the vector field given by Definition \ref{Ondelette_base}.   For  $1\leq p \leq +\infty$ there exist two constants $0<c_{p} < {c'}_{p}$, which only depend on   $p$ and $\vec{\phi}$, such that, for all sequence  $\lambda=(\lambda_k)_{k\in \mathbb{Z}^3}\in \ell^{p}(\mathbb{Z}^3)$ we have the \emph{almost-orthogonality property}:
\begin{equation*}
c_{p}\Vert \lambda\Vert_{\ell^{p}(\Zt)} \leq \left\Vert \sum_{k\in \mathbb{Z}^3}\lambda_k  \vec{\phi}(\cdot -k)\right\Vert_{L^p(\Rt)}\leq {c'}_{p}\Vert \lambda\Vert_{\ell^{p}(\Zt)}.
\end{equation*}
\end{Lemme}
For a proof of this fact and for more references about wavelets' properties see the books \cite{KahLer} (Part II, Chapter $6$)  or  \cite{Meyer}  (Chapter III).\\

Once this function $\vec{\phi}$ is fixed, we will construct the force $\fe$. Let $\theta\ell_0>0$ be an energy input scale which is fixed from now on. Let $ \ell \geq \theta\ell_0$ be a parameter from which we consider the cube $[-\ell,\ell]^3\subset \Rt$. As the energy input scale $\ell_0$ is such that $\frac{\ell}{\theta\ell_0}\geq 1$  then, in the cube $[-\ell,\ell]^3$, we shall consider all the points of the form $\theta\ell_0k$ where $\vert \theta\ell_0k\vert\leq \ell$ with $k\in \mathbb{Z}^{3}$.  Thus following an idea of \cite{DoerFoias} (which was given in the periodic setting), we will construct the force $\fe$ by translations of the function $\vec{\phi}$ to each point $\theta\ell_0 k$  and by dilatation to the scale $\frac{1}{\theta\ell_0}$.  Hence we get the formula 
\begin{equation}\label{Def_Force}
\fe(x)=A\sum_{\vert \ell_0k\vert \leq \ell}\vec{\phi}\left(\frac{x-\theta\ell_0k}{\theta\ell_0}\right),
\end{equation} where the parameter $A>0$ is the amplitude of the external force.
\begin{Lemme}\label{Prop-force}  Let  $\fe$ be the force given in (\ref{Def_Force}) above.
\begin{enumerate}
\item[1)]  We have $supp\,\left(\fe\right)\subset \lbrace \xi \in \Rt: \frac{1}{ 10 \theta\ell_0} \leq \vert \xi \vert\leq  \frac{1}{\theta\ell_0}\rbrace$. 
\item[2)] For all  $1\leq p \leq +\infty$ there exist two constants $0<C_{p}< {C'}_{p}$ such that
\begin{equation*}
C_{p}\, A \ell^{\frac{3}{p}} \leq \| \fe \|_{L^p}\leq {C'}_{p} \,A \ell^{\frac{3}{p}}.
\end{equation*}
\end{enumerate} 
\end{Lemme}
\pv  
\begin{enumerate}
\item[1)] Due to point $2)$ of Definition \ref{Ondelette_base} and formula (\ref{Def_Force}) the proof of this frequency localization is  straightforward. 
\item[2)]  Consider the  sequence  $(\lambda_k)_{k\in \Zt}$  defined as
\begin{equation*}
\lambda_k=\left\lbrace
\begin{array}{rl} 1&\text{if}\quad \vert \theta\ell_0 k \vert\leq \ell, \\[2mm]
0&\text{otherwise.}\end{array}\right.
\end{equation*}  
Then we can write 
\begin{equation*}
\fe(x)=A \sum_{\vert \theta\ell_0k\vert\leq \ell}  \vec{\phi}\left(\frac{x-\theta\ell_0k}{\theta\ell_0}\right)=A\sum_{k\in \Zt}\lambda_k  \vec{\phi} \left(\frac{x}{\theta\ell_0}- k\right).
\end{equation*}
Thus, for $1\leq p \leq+\infty$, taking the $L^p$ norm in the identity above and by  Lemma \ref{prop_ondelette_base}  there exist two constants  $c_{p},{c'}_{p}>0$ such that
\begin{equation}\label{est_tech_1}
c_{p}\,A (\theta \ell_0)^{\frac{3}{p}} \Vert \lambda \Vert_{\ell^{p}} \leq A (\theta \ell_0)^{\frac{3}{p}}\left\Vert \sum_{k\in \mathbb{Z}^3}\lambda_k  \vec{\phi}(\cdot -k)\right\Vert_{L^p}\leq {c'}_{p}\,A (\theta \ell_0)^{\frac{3}{p}} \Vert \lambda \Vert_{\ell^{p}},
\end{equation}
On the other hand, by definition of the sequence   $(\lambda_k)_{k\in \Zt}$ we may see that $\Vert \lambda\Vert^{p}_{L^p}$ is the number of points of the form $\theta\ell_0k$ in the cube $[-\ell,\ell]^3$ and thus, there exist   two constants $c,c'>0$ (always independent of the  parameters $\ell_0$, $\ell$, $\theta$ and $A$ ) such that  
\begin{equation*}
c\left(\frac{\ell}{\theta \ell_0}\right)^{\frac{3}{p}} \leq \Vert \lambda\Vert_{\ell^{p}}\leq c'\left(\frac{\ell}{\theta \ell_0}\right)^{\frac{3}{p}},
\end{equation*}
and thus by expression (\ref{est_tech_1})  we obtain the desired estimate. \finpv
\end{enumerate}
This lemma gives us all the properties required for the external force $\fe$. Indeed, in point $1)$ we may see that this functions is localized at the frequencies $\frac{1}{10 \theta\ell_0} \leq \vert \xi \vert \leq \frac{1}{\theta\ell_0}$, and moreover, boy point $2)$ and we shall  write from now on
\begin{equation}\label{Formula_EstimateForce}
\| \fe \|_{L^p} \simeq A \ell^{\frac{3}{p}},
\end{equation} 
hence we get the following estimations:  if $p=+\infty$ then we obtain $\Vert \fe \Vert_{L^{\infty}}\simeq A,$ and we may see that  the amplitude of the external force $\fe$ is indeed given by the parameter $A$. Thus, we set $A= \frac{\nu^2}{\ell^{3}_{0}}$ and by expression (\ref{Def_Force}) we write, from now on 
\begin{equation}\label{fe}
\fe(x)=\frac{\nu^2}{\ell^{3}_{0}}\, \sum_{\vert\theta \ell_0k\vert \leq \ell}\vec{\phi}\left(\frac{x-\theta\ell_0k}{\theta\ell_0}\right).
\end{equation}
Moreover, if $p=2$ the we obtain $\Vert \fe \Vert_{L^2} \simeq  \frac{\nu^2}{\ell^{3}_{0}}\,  \ell^{\frac{3}{2}}$, where we may that the quantity $\Vert \fe \Vert_{L^2}$ is driven only by means of the parameter $\ell$ and the property $\Vert \fe \Vert_{L^2} \gg 1$ is equivalent to  $\ell \gg 1$.\\
	 
With this external force we can reach large values for the Grashof number $Gr$. Indeed, recall first that the parameter $\gamma>0$ is given in (\ref{gamma}) by $\gamma= \frac{\Vert \fe \Vert_{L^{\infty}}}{c_0 \left(\frac{1}{\ell_0}\right)^{\frac{3}{2}} \Vert \fe \Vert_{L^2}}$ and by since we have  $\Vert \fe \Vert_{L^{\infty}} \simeq \frac{\nu^2}{\ell^{3}_{0}}$ and $\Vert \fe \Vert_{L^2} \simeq  \frac{\nu^2}{\ell^{3}_{0}}\,  \ell^{\frac{3}{2}}$, then we get $\gamma \simeq \frac{\ell^{\frac{3}{2}}_{0}}{\ell^{\frac{3}{2}}}$. But, since the number $G_0$ is fixed then by (\ref{estimate_Grashof_G0}) we have $Gr \simeq \frac{1}{\gamma^4}$ and thus we get $Gr \simeq \frac{\ell^6}{\ell^{6}_{0}}$, where we may see that $Gr\gg 1$ if $\ell\gg 1$. \\

However, even if the Grashof number can be large enough,  we will see   the Taylor scale $\ell_T$ verifies $\ell_T \simeq \ell_0$, as stated in Theorem \ref{Th-non-turb},  and then the turbulence is annihilated. Indeed, recall that the Taylor scale $\ell_T$ is defined by $\ell_T= \left( \frac{\nu U^2}{\varepsilon}\right)^{\frac{1}{2}}$. In order to estimate $U$ remark that by Proposition \ref{lemme_estimate_U_F_L} we have $ \frac{U^2}{L} \simeq F$ hence we write $U \simeq (L F)^{\frac{1}{2}}$. But  we have $L= \frac{\ell_0}{\gamma}$ and $F= \frac{\Vert \fe \Vert_{L^2}}{\ell^{\frac{3}{2}}_{0}}$ (see (\ref{LC}) and (\ref{F})) and by estimates on $\gamma$ and $\Vert \fe \Vert_{L^2}$ above  we have $L \simeq \frac{\ell^{\frac{3}{2}}}{\ell^{\frac{1}{2}}_{0}}$ and $F \simeq \frac{\nu^2 \ell^{\frac{3}{2}}}{\ell^{\frac{9}{2}}_{0}}$ respectively. Then we get  $U \simeq  \frac{\nu \ell^{\frac{3}{2}}}{\ell^{\frac{5}{2}}_{0}}$. In order to estimate $\varepsilon$ remark that by Theorem \ref{Main_Result} we have $\varepsilon \simeq \frac{U^3}{L}$, hence we get  $\varepsilon \simeq \frac{\nu^3 \ell^{3}}{\ell^{7}_{0}}$.  With these estimates on $U$ and $\varepsilon$ at hand, we get back to the expression of the Taylor scale and we write 
$$ \ell_T = \left(\frac{\nu U^2}{\varepsilon}\right)^{\frac{1}{2}} \simeq \left( \nu  \frac{ \nu^2  \ell^{3}}{\ell^{5}_{0}} \frac{\ell^{7}_{0}}{ \nu^3 \ell^{3}}\right)^{\frac{1}{2}}=\ell_0.$$ 
%%%%%%%%%%%%%%%%%%%%%%%%%%%%%%%%%%%%%%%
\begin{appendices}
\section{Appendix}\label{AppendixA}
Let $\vu \in L^{\infty}_{loc}([0,+\infty[, L^{2}(\Rt)) \cap L^{2}_{loc}([0,+\infty[, \dot{H}^{1}(\Rt))$ be a Leray's weak solution of the classical Navier-Stokes equations. We give here a short proof of the fact that the energy dissipation rate $\varepsilon$ defined by means of $\vu$ in (\ref{varepsilon,U}) verifies $\varepsilon<+\infty$. As mentioned in the introduction, this fact relies on the classical energy inequality: 
$$ \Vert \vu(t,\cdot)\Vert^{2}_{L^2}+2\nu\int_{0}^{t} \Vert 
\vu(s,\cdot)\Vert^{2}_{\dot{H}^{1}}ds\leq \Vert \vu_0 \Vert^{2}_{L^2}+2\int_{0}^{t} \langle \fe, \vu(s,\cdot) \rangle_{\dot{H}^{-1}\times \dot{H}^1} ds.$$ Indeed,  since  $\Vert \vu(t,\cdot)\Vert^{2}_{L^2}$ is a positive quantity  we can write 
$$ 2\nu\int_{0}^{t} \Vert 
\vu(s,\cdot)\Vert^{2}_{\dot{H}^{1}}ds\leq \Vert \vu_0 \Vert^{2}_{L^2}+2\int_{0}^{t} \langle \fe, \vu(s,\cdot) \rangle_{\dot{H}^{-1}\times \dot{H}^1} ds.$$ Then, remark that by the hypothesis on $\vu$ and $\fe$,  the second term in the right-hand side can be  estimated  as follows: 
\begin{eqnarray}\label{estim1} \nonumber 
2\int_{0}^{t} \langle \fe, \vu(s,\cdot) \rangle_{\dot{H}^{-1}\times \dot{H}^1} ds  &\leq &2 \int_{0}^{t} \Vert \fe \Vert_{\dot{H}^{-1}} \Vert \vu(s,\cdot)\Vert_{\dot{H}^{1}} ds \leq \int_{0}^{t} (\frac{\Vert \fe \Vert^{2}_{\dot{H}^{-1}}}{\nu} + \nu \Vert \vu(s,\cdot)\Vert^{2}_{\dot{H}^{1}}) ds \\
& \leq & t \frac{\Vert \fe \Vert^{2}_{\dot{H}^{-1}}}{\nu} + \nu \int_{0}^{t} \Vert \vu (s,\cdot)\Vert^{2}_{\dot{H}^{1}}ds, 
\end{eqnarray} and getting back to the previous estimate on the quantity $\ds{2\nu\int_{0}^{t} \Vert 
\vu(s,\cdot)\Vert^{2}_{\dot{H}^{1}}ds}$, we obtain 
$$ 2\nu\int_{0}^{t} \Vert 
\vu(s,\cdot)\Vert^{2}_{\dot{H}^{1}}ds\leq \Vert \vu_0 \Vert^{2}_{L^2} \leq   t \frac{ \Vert \fe \Vert^{2}_{\dot{H}^{-1}}}{\nu} + \nu \int_{0}^{t} \Vert \vu (s,\cdot)\Vert^{2}_{\dot{H}^{1}}ds,$$ hence we write 
\begin{equation}\label{estim2}
 \nu\int_{0}^{t} \Vert  \vu(s,\cdot)\Vert^{2}_{\dot{H}^{1}}ds \leq \Vert \vu_0 \Vert^{2}_{L^2} +   t \frac{\Vert \fe \Vert^{2}_{\dot{H}^{-1}}}{\nu}.
\end{equation} Finally, we divide each term by $t$ and  $\ell^{3}_{0}$, then we take the limit $\ds{\limsup_{t \longrightarrow+\infty}}$, and by definition of the quantity $\varepsilon$ we have 
$$ \varepsilon =\nu\limsup_{t\longrightarrow +\infty} \frac{1}{t}\int_{0}^{t}\Vert \vu(s,\cdot)\Vert^{2}_{\dot{H}^1}\frac{ds}{\ell^3_{0}} \leq \frac{\Vert \fe \Vert^{2}_{\dot{H}^{-1}}}{\nu \ell^{3}_{0}}<+\infty.$$
%%%%%%%%%%%%%%%%%%%%%%%%%%%%%%%%%%%%%
\section{Appendix}\label{AppendixB}
 For  $L>0$  consider the cube $[0,L]^{3}\subset \Rt$,  and let  $\vu \in L^{\infty}_{loc}([0,+\infty[, L^{2}([0,L]^3)) \cap L^{2}_{loc}([0,+\infty[, \dot{H}^{1}([0,L]^3))$ be a Leray's weak solution of the classical and periodic Navier-Stokes equations. We will prove that in the periodic setting the characteristic velocity defined in (\ref{varepsilon,U}) verifies $U<+\infty$. \\
 
Start with the estimate given in (\ref{estim2}), which is also verified in the periodic setting, and we study now the term $ \ds{ \nu\int_{0}^{t} \Vert  \vu(s,\cdot)\Vert^{2}_{\dot{H}^{1}}ds}$. Indeed, by the following Poincaré inequality $\ds{\Vert \vu(t,\cdot) \Vert_{L^2} \leq \frac{ L}{2\pi} \Vert \vec{\nabla} \otimes \vu(t,\cdot)\Vert_{L^2}}$, we can write $\ds{\Vert \vu(t,\cdot) \Vert^{2}_{L^2} \leq \frac{ L^2}{4\pi^{2}} \Vert \vec{\nabla} \otimes \vu(t,\cdot)\Vert^{2}_{L^2} \leq c \frac{ L^{2}}{4\pi^{2}} \Vert \vu(s,\cdot)\Vert^{2}_{\dot{H}^{1}}}$, hence  we obtain 
$$ \int_{0}^{t} \Vert \vu(s,\cdot) \Vert^{2}_{L^2} ds \leq c  \frac{ L^{2}}{4\pi^{2}} \int_{0}^{t} \Vert \vu(s,\cdot) \Vert^{2}_{\dot{H}^{1}} ds,$$ and thus we have the estimate from below $\ds{ \frac{\nu }{c} \frac{4 \pi^2}{L^2} \int_{0}^{t} \Vert \vu(s,\cdot)\Vert^{2}_{L^2} ds \leq   \nu\int_{0}^{t} \Vert  \vu(s,\cdot)\Vert^{2}_{\dot{H}^{1}}ds}$. Then, by estimate  (\ref{estim2}) we can write 
$$ \frac{\nu }{c} \frac{4 \pi^2}{L^2} \int_{0}^{t} \Vert \vu(s,\cdot)\Vert^{2}_{L^2} ds \leq   \Vert \vu_0 \Vert^{2}_{L^2} +   t \frac{\Vert \fe \Vert^{2}_{\dot{H}^{-1}}}{\nu},$$ where, diving by $t$ and $\ell^{3}_{0}$, and taking the limit $\ds{\limsup_{t \longrightarrow+\infty}}$, by definition of the quantity $U$ we have $$ \frac{\nu }{c} \frac{4 \pi^2}{L^2} U^2 \leq  \frac{\Vert \fe \Vert^{2}_{\dot{H}^{-1}}}{\nu \ell^{3}_{0}},$$ hence we have $U<+\infty$.

%%%%%%%%%%%%%%%%%%%%%%%%%%%%%%%%%%%%% 
\section{Appendix}\label{AppendixC}
Remark that the control in time given in (\ref{control1}) follows directly from the energy inequality (\ref{energy-classical}) and the inequality given in (\ref{estim1}). Indeed, by these inequalities we can write 
$$ \Vert \vu(t,\cdot)\Vert^{2}_{L^2}+2\nu\int_{0}^{t} \Vert 
\vu(s,\cdot)\Vert^{2}_{\dot{H}^{1}}ds\leq \Vert \vu_0 \Vert^{2}_{L^2}+t \frac{\Vert \fe \Vert^{2}_{\dot{H}^{-1}}}{\nu} + \nu \int_{0}^{t} \Vert \vu (s,\cdot)\Vert^{2}_{\dot{H}^{1}}ds,$$ hence we obtain 
$$ \Vert \vu(t,\cdot)\Vert^{2}_{L^2}+\nu\int_{0}^{t} \Vert 
\vu(s,\cdot)\Vert^{2}_{\dot{H}^{1}}ds\leq \Vert \vu_0 \Vert^{2}_{L^2}+t \frac{\Vert \fe \Vert^{2}_{\dot{H}^{-1}}}{\nu},$$ but $\ds{\nu\int_{0}^{t} \Vert 
\vu(s,\cdot)\Vert^{2}_{\dot{H}^{1}}ds}$ is a positive quantity and then we have 
$$  \Vert \vu(t,\cdot)\Vert^{2}_{L^2} \leq \Vert \vu_0 \Vert^{2}_{L^2}+t \frac{\Vert \fe \Vert^{2}_{\dot{H}^{-1}}}{\nu},$$ which is the control in time given in (\ref{control1}).
%%%%%%%%%%%%%%%%%%%%%%%%%%%%%%%%%%%%%
\end{appendices}

%%%%%%%%%%%%%%%%%%%%%%%%%%%%%%%%%%%%%%
%%%%%%%%%%%%%%%%%%%%%%%%%%%%%%%%%%%%%%


\begin{thebibliography}{20}
%%%%%%%%%%%%%%%%%%%%%%%%%%%%%%%%%%%%%%
\bibitem{Chesk} A. Cheskidov, S. Friedlander \&  N. Pavlovi\'c.  \emph{An inviscid dyadic model of turbulence: the global attractor}. AIMS.,    26: 781 - 794, (2010).    
\bibitem{Const2} P. Constantin, Q. Nie \& S. Tanveer. \emph{Bounds for second order structure functions and energy spectrum in turbulence}.  Physics of fluids., 11, (1999).
\bibitem{Const} P. Constantin. \emph{Euler equations Navier-Stokes equations and turbulence, Mathematical Foundation of Turbulent Viscous Flows}, Vol. 1871 of the series Lecture Notes in Mathematics: 1-43, (2005).
\bibitem{DasGru} R. Dascaliuc and Z. Gruli\'c. \emph{On energy cascades in the forced 3D Navier-Stokes equations}. J Nonlinear Sci.,  22:683-715,  (2016). 
\bibitem{DoerFoias} C. Doering and C. Foias. \emph{Energy dissipation in body-forced turbulence}, Journal of Fluid Mechanics., 467:289- 306, (2002).
\bibitem{FMRT1} C. Foias, O. Manley, R. Rosa \& R. Temam. \emph{Estimates for the energy cascade in three-dimensional turbulent flows}.  Comptes Rendus Acad. Sci. Paris, Série I., 332:509-514, (2001).
\bibitem{FMRT2} C. Foias, M. S. Jolly, O. Manley, R. Rosa \& R. Temam.  \emph{Kolmogorov theory via finite-time averages}. Physica D: Nonlinear Phenomena Volume., 212:245-270  (2005). 
\bibitem{Houel} G. Hou\"el. \emph{Estimates on the energy dissipation rate in kinetic turbulence}. Research project in the laboratory  MEC559-MEC569, Ecole Polytechnique, (2004).
\bibitem{JacTab} L. Jacquin and P. Tabeling.  \emph{Turbulence et Tourbillons}.  Majeure 1, MEC555, lectures notes  Ecole Polytechnique,  (2006). 
\bibitem{OJarrin} O. Jarr\'in.  \emph{Some remarks on the Kolmogorov dissipation law in the deterministic framework of Navier-Stokes equations in $\Rt$}. In process. 
\bibitem{OscarTesis} O. Jarr\'in. \emph{Deterministic descriptions of turbulence in the Navier-Stokes equations}. Ph.D. thesis at Universit\'e Paris-Saclay, (2018). hal-01821762.
\bibitem{KahLer} J.P. Kahane and P.G.  Lemarié-Rieusset. \emph{S\'erie de Fourier et ondelettes}. Cassini, Paris, (1998).
\bibitem{Kolm1} A. N. Kolmogorov.  \emph{The local structure of turbulence in incompressible viscous fluid for very large Reynolds number}. Dokl. Akad. Nauk SSSR, 30 (1941).
\bibitem{Kolm2} A. N. Kolmogorov.  \emph{On Degeneration of Isotropic turbulence in an incompressible viscous liquid}. Dokl. Akad. Nauk SSSR 31 (1941).
\bibitem{PGLR2} P.G. Lemarié-Rieusset. \emph{Recent developments in the Navier-Stokes problem}, Chapman \& Hall/CRC, (2002).
\bibitem{PGLR1} P.G. Lemarié-Rieusset. \emph{The Navier-Stokes Problem in the 21st Century}, Chapman \& Hall/CRC, (2016).
\bibitem{Meyer} Y. Meyer. \emph{Ondelettes et Operateurs}, Tome 1.  Hermann, (1990). 
\bibitem{Obuk} A.M. Obukhoff.  \emph{On the energy distribution in the spectrum of a turbulent flow}. Dokl. Akad. Nauk SSSR, 32, (1941).
\bibitem{OttoRamos} F. Otto and F. Ramos. \emph{Universal bounds for the Littlewood-Paley first-order moments of the 3-D Navier-Stokes equations}, 2009. Commun. Math. Phys., 300: 301-315, (2010).
\bibitem{Tennekes}  H. Tennekes and  J.L. Lumley. \emph{A First Course in Turbulence}. The MIT Press (1972).
\bibitem{Vigneron} F. Vigneron. \emph{Free turbulence on $\mathbb{R}^3$ and $\mathbb{T}^3$}. Dynamics of PDE., 7:107-160,  (2010).
\bibitem{Wilcox} D. Wilcox.  \emph{Turbulence model for CFD}. DCW Industries, Inc (1994).
\end{thebibliography}
\end{document}